\documentclass[11pt]{amsart}

\usepackage{amsmath, amsthm, amssymb, amsfonts, mathrsfs, amscd}

\def\CC{{\mathbb C}}

\def\KK{{\mathbb K}}

\def\PP{{\mathbb P}}

\def\RR{{\mathbb R}}

\def\hhat{{\hat h}}
\def\Hhat{{\hat H}}

\def\0{{\mathbf 0}}
\def\1{{\mathbf 1}}

\def\e{{\mathbf e}}

\def\x{{\mathbf x}}
\def\y{{\mathbf y}}

\def\Acal{{\mathcal A}}
\def\Bcal{{\mathcal B}}

\def\Ecal{{\mathcal E}}

\def\Lcal{{\mathcal L}}
\def\Mcal{{\mathcal M}}
\def\Ocal{{\mathcal O}}

\def\Scal{{\mathcal S}}

\def\hhat{{\widehat h}}

\def\Kbar{{\bar K}}
\def\kbar{{\bar k}}

\def\End{\mathrm{End}}
\def\Fix{\mathrm{Fix}}
\def\GL{\mathrm{GL}}

\def\Id{\mathrm{Id}}

\def\ord{\mathrm{ord}}
\def\PGL{\mathrm{PGL}}
\def\Per{\mathrm{Per}}
\def\Perm{\mathrm{Perm}}
\def\Pic{\mathrm{Pic}}
\def\Preper{\mathrm{PrePer}}
\def\Res{\mathrm{Res}}

\def\Spec{\mathrm{Spec}}
\def\St{\mathrm{St}}
\def\Sym{\mathrm{Sym}}

\newcommand{\rst}[1]{\ensuremath{{\mathbin\upharpoonright}\raise-.5ex\hbox{$#1$}}}

\theoremstyle{plain}

\newtheorem{thm}{Theorem}
\newtheorem{cor}[thm]{Corollary}
\newtheorem{prop}[thm]{Proposition}
\newtheorem{lem}[thm]{Lemma}

\begin{document}

\title[Isotriviality and potential good reduction]{Isotriviality is equivalent to potential good reduction for endomorphisms of $\PP^N$ over function fields}

\dedicatory{Dedicated to Paul Roberts on the occasion of his 60th birthday.}
\date{June 8, 2008}

\author{Clayton Petsche}
\email{cpetsche@gc.cuny.edu}
\address{Clayton Petsche; Ph.D. Program in Mathematics; CUNY Graduate Center; 365 Fifth Avenue; New York, NY 10016-4309 U.S.A.}

\author{Lucien Szpiro}
\email{lszpiro@gc.cuny.edu}
\address{Lucien Szpiro; Ph.D. Program in Mathematics; CUNY Graduate Center; 365 Fifth Avenue; New York, NY 10016-4309 U.S.A.}

\author{Michael Tepper}
\email{michael.tepper@gmail.com}
\address{Michael Tepper; Ph.D. Program in Mathematics; CUNY Graduate Center; 365 Fifth Avenue; New York, NY 10016-4309 U.S.A.}

\subjclass[2000]{14G99,  14H05} 

\begin{abstract}
Let $K=k(C)$ be the function field of a complete nonsingular curve $C$ over an arbitrary field $k$.  The main result of this paper states that a morphism $\varphi:\PP^N_K\to\PP^N_K$ is isotrivial if and only if it has potential good reduction at all places $v$ of $K$; this generalizes results of Benedetto for polynomial maps on $\PP^1_K$ and Baker for arbitrary rational maps on $\PP^1_K$.  We offer two proofs: the first uses algebraic geometry and geometric invariant theory, and it is new even in the case $N=1$.  The second proof uses non-archimedean analysis and dynamics, and it more directly generalizes the proofs of Benedetto and Baker.  We will also give two applications.  The first states that an endomorphism of $\PP^N_K$ of degree at least two is isotrivial if and only if it has an isotrivial iterate.  The second gives a dynamical criterion for whether (after base change) a locally free coherent sheaf $\Ecal$ of rank $N+1$ on $C$ decomposes as a direct sum $\Lcal\oplus\dots\oplus\Lcal$ of $N+1$ copies of the same invertible sheaf $\Lcal$.
\end{abstract}

\maketitle


\begin{section}{Introduction}

Let $K=k(C)$ be the function field of a complete nonsingular curve $C$ over an arbitrary field $k$, and let $\varphi:\PP^N_K\to\PP^N_K$ be a morphism.  We say $\varphi$ is {\em trivial} if it is defined over the constant field $k$.  More generally, we say $\varphi$ is {\em isotrivial} if there exists a finite extension $K'/K$ such that the induced morphism $\PP^N_{K'}\to\PP^N_{K'}$ is defined over the algebraic closure $k'$ of $k$ in $K'$.

In the study of dynamical systems arising from the iteration of morphisms $\varphi:\PP^N_K\to\PP^N_K$ of degree at least two, the isotrivial and non-isotrivial cases exhibit very different behavior.  For example, in the one-dimensional case Baker \cite{Baker} showed that if $\varphi:\PP^1_K\to\PP^1_K$ is a non-isotrivial morphism of degree at least two, then $\PP^1(K)$ has only finitely many small points with respect to the Call-Silverman canonical height function $\hhat_\varphi$ associated to $\varphi$.  In particular, it follows that $\PP^1(K)$ contains only finitely many $\varphi$-preperiodic points, and that a point $P\in\PP^1(\Kbar)$ is preperiodic if and only if $\hhat_\varphi(P)=0$.  These results had been previously established in the special case of polynomial endomorphisms of $\PP^1_K$ -- that is, endomorphisms with a totally ramified fixed point -- by Benedetto \cite{Benedetto}.  Using entirely different techniques from model theory, Chatzidakis-Hrushovski \cite{ChatzidakisHrushovski} have recently generalized Baker's results to endomorphisms of $\PP^N_K$.

The situation is quite different for isotrivial endomorphisms.  Suppose for example that the constant field $k$ is algebraically closed and that $\varphi:\PP^1_K\to\PP^1_K$ is a morphism defined over $k$ with $\deg(\varphi)\geq2$; then all of the (infinitely many) $\varphi$-preperiodic points in $\PP^1(\Kbar)$ are defined over $k$, and so they are $K$-rational.  Moreover, since $\varphi$ is isotrivial the canonical height $\hhat_\varphi$ coincides with the naive height, so unless $k$ is an algebraic closure of a finite field, $\PP^1(K)$ may contain non-preperiodic points having canonical height zero.

A key ingredient in the results of Baker and Benedetto on non-isotrivial endomorphisms of $\PP^1_{K}$ is a characterization of isotriviality in terms of purely local conditions; see Baker \cite{Baker} Thm. 1.9, and also the related Prop. 6.1 of Benedetto \cite{Benedetto}.  The main result of this paper generalizes this criterion to endomorphisms of $\PP^N_{K}$.  In order to state the theorem, we first recall that the set $M_K$ of places of $K$ can be naturally identified with the set of closed points on the curve $C$.  Given a place $v\in M_K$, denote by $\Ocal_v\subset K$ the ring of regular functions at $v$.  We say $\varphi$ has {\em good reduction} at $v$ if there exists a choice of homogeneous coordinates $(x_0,x_1,\dots,x_N)$ on $\PP^N_{K}$ such that $\varphi$ extends to an endomorphism of the associated integral model $\PP^N_{\Ocal_v}$ of $\PP^N_{K}$.  We say $\varphi$ has {\em potential good reduction} at $v$ if there exists a finite extension $K'/K$ and a place $v'$ of $K'$ over $v$ such that $\varphi$ has good reduction at $v'$.

\begin{thm}\label{MainTheorem}
Let $K=k(C)$ be a function field, and let $\varphi:\PP^N_K\to\PP^N_K$ be a morphism.  Then $\varphi$ is isotrivial if and only if $\varphi$ has potential good reduction at all places $v$ of $K$.
\end{thm}

We first remark that the ``only if'' direction of the theorem is easy; see $\S$~\ref{OnlyIfSect}.  When $\deg(\varphi)=0$ there is nothing to prove, since all constant morphisms are trivial and have good reduction at all places.  The $\deg(\varphi)=1$ case is a simple exercise in linear algebra; we will give the proof in $\S$~\ref{DegOneCaseSect}.  Thus, the interesting part of Theorem~\ref{MainTheorem} is the ``if'' direction when $\deg(\varphi)\geq2$.  

As a reflection of the variety of techniques which are commonly used to study algebraic dynamics over global fields, we will give two very different proofs of Theorem~\ref{MainTheorem}.  Our first proof uses algebraic geometry and standard facts from geometric invariant theory, and it is new even in the one-dimensional case.  Given integers $N\geq1$ and $d\geq2$, we first show that the space $\Mcal_{N,d}$ parametrizing morphisms $\varphi:\PP^N_k\to\PP^N_k$ with $\varphi^*\Ocal(1)\simeq\Ocal(d)$ exists as an affine $k$-variety (when $N=1$ this follows immediately from Silverman \cite{SilvermanSpace}, Thm. 1.1).  A morphism $\varphi:\PP^N_K\to\PP^N_K$ over the function field $K=k(C)$ with everywhere potential good reduction induces a regular map $C\to\Mcal_{N,d}$, which must be constant since $\Mcal_{N,d}$ is affine.  It follows that $\varphi$ is isotrivial.

Our second proof of Theorem~\ref{MainTheorem} uses non-archimedean analysis and dynamics, and it more directly generalizes the proofs given by Benedetto \cite{Benedetto} and Baker \cite{Baker}.  We consider the local homogeneous filled Julia set $F_{\Phi,v}$ associated to each place $v$ of $K$ and each model $\Phi$ for $\varphi$; this is a certain dynamical invariant of $\Phi$ which detects good reduction at $v$.  The key step is to define a notion of homogeneous transfinite diameter in order to measure the size of the set $F_{\Phi,v}$, and to show that these numbers satisfy a product formula over all places $v$ of $K$.  Selecting a globally defined model $\Phi$ for $\varphi$ with certain favorable properties, we use the hypothesis of everywhere potential good reduction along with the product formula to show that $\Phi$ must be defined over the constant field $k$.

Despite the different techniques used in our two approaches, they nevertheless share several ingredients in common.  For example, both proofs use basic facts about the resultant of homogeneous maps (which we will review in $\S$~\ref{ResSect}), and both proofs use the fact that the curve $C$ has no non-constant regular functions.  Moreover, both proofs make essential use of a basic result in algebraic dynamics on the Zariski-density of preperiodic points (we will discuss this further in $\S$~\ref{AlgDynSys}).

In $\S$~\ref{ApplicationsSect} we will give two applications of Theorem~\ref{MainTheorem}.  The first is a result stating that an endomorphism of $\PP^N_K$ of degree at least two is isotrivial if and only if it has an isotrivial iterate.  The second application gives a dynamical criterion for whether (after base change) a locally free coherent sheaf $\Ecal$ of rank $N+1$ on $C$ decomposes as a direct sum $\Lcal\oplus\dots\oplus\Lcal$ of $N+1$ copies of the same invertible sheaf $\Lcal$ on $C$.  When $k=\CC$, Amerik \cite{Amerik} has obtained a similar result with the curve $C$ replaced by a smooth projective base $B$ of arbitrary dimension.

It would be interesting to investigate to what extent Theorem~\ref{MainTheorem} can be extended to more general polarized algebraic dynamical systems (in the sense of $\S$~\ref{AlgDynSys}).  We caution, however, that a naive restatement of Theorem~\ref{MainTheorem} in this general setting is false, as there do exist non-isotrivial endomorphisms $\varphi:X\to X$ of projective $K$-varieties $X$ with good reduction at every place $v$ of $K$, in the sense that $\varphi$ extends to a projective endomorphism $\varphi_v:X_{\Ocal_v}\to X_{\Ocal_v}$ of an integral model $X_{\Ocal_v}$ for $X$.  For example, let $k$ be an algebraically closed field of characteristic zero, let $g$ and $n$ be large  positive integers, and let $\Acal_{g,n}$ denote the (fine) moduli space of principally polarized abelian varieties over $k$ of dimension $g$ with level-$n$ structure.  It is well-known that $\Acal_{g,n}$ contains a complete nonsingular curve $C$.  (This follows from the fact that the boundary $\overline{\Acal}_{g,n}\setminus\Acal_{g,n}$ of $\Acal_{g,n}$ inside of its Satake compactification $\overline{\Acal}_{g,n}$ has codimension strictly greater than one when $g$ is large enough; see for example Cartan \cite{Cartan} or Kodaira \cite{Kodaira}.)  The resulting abelian scheme $A\to C$ has as its generic fiber an abelian variety $A_K$ over $K=k(C)$, and the doubling map $[2]:A_K\to A_K$ gives rise to a non-isotrivial polarized algebraic dynamical system over $K$ with everywhere good reduction.

{\em Acknowledgments.}  The authors would like to thank Rob Benedetto, Antoine Chambert-Loir, Charles Favre, and Felipe Voloch for their helpful suggestions.

\end{section}


\begin{section}{Preliminaries}\label{PrelimSect}

\begin{subsection}{Review of polarized algebraic dynamical systems}\label{AlgDynSys}

Let $k$ be a field.  A polarized algebraic dynamical system over $k$ is triple $(X,\varphi,\Lcal)$, where $X$ is a projective $k$-variety, $\varphi:X\to X$ is a morphism, and $\Lcal$ is an ample invertible sheaf on $X$ such that $\varphi^*\Lcal\simeq\Lcal^{\otimes d}$ for some $d\geq2$.  We will now recall several standard definitions and facts about polarized algebraic dynamical systems; for more background see the surveys by Zhang \cite{Zhang} and Chambert-Loir \cite{ChambertLoir}.

Given a point $x\in X(\bar{k})$, we say $x$ is {\em fixed} if $\varphi(x)=x$, {\em periodic} if $\varphi^n(x)=x$ for some $n\geq1$, and {\em preperiodic} if $\varphi^m(x)$ is periodic for some $m\geq0$.  Let $\Fix(\varphi)$, $\Per(\varphi)$, and $\Preper(\varphi)$ denote the sets of fixed, periodic, and preperiodic points in $X(\bar{k})$, respectively.  Given integers $n\geq1$ and $m\geq0$ we denote by
\begin{equation*}
\begin{split}
\Per_n(\varphi) & = \{x\in X(\bar{k})|\varphi^{n}(x)=x\} \\
\Preper_{n,m}(\varphi) & = \{x\in X(\bar{k})|\varphi^{n+m}(x)=\varphi^m(x)\}
\end{split}
\end{equation*}
the sets of periodic points of period $n$, and preperiodic points of type $(n,m)$, respectively.

\begin{prop}\label{DegreeFormula}
Let $(X,\varphi,\Lcal)$ be a polarized algebraic dynamical system, and assume that $X$ is geometrically integral.  Then:
\begin{quote}
{\bf (a)}  The morphism $\varphi$ is finite and $\deg(\varphi)=d^{\dim(X)}$. \\
{\bf (b)}  The sets $\Fix(\varphi)$, $\Per_n(\varphi)$, and $\Preper_{n,m}(\varphi)$ are finite.  \\
{\bf (c)}  The set $\Preper(\varphi)$ is Zariski-dense in $X$.
\end{quote}
\end{prop}

\begin{proof}
{\bf (a)}  Suppose on the contrary that $\varphi$ is not finite.  Since a projective morphism is finite if and only if it has finite fibers, this means that there exists a point $x\in X(\kbar)$ such that $\varphi^{-1}(x)$ contains an irreducible curve $Z$.  Pushing forward the intersection product of $Z$ with the first Chern class $c_1(\varphi^*\Lcal)$ of $\varphi^*\Lcal$, we have
\begin{equation}\label{ProjForm}
\varphi_*(Z.c_1(\varphi^*\Lcal))=\varphi_*(Z).c_1(\Lcal)
\end{equation} 
by the projection formula.  We have $Z.c_1(\varphi^*\Lcal)=Z.c_1(\Lcal^{\otimes d})=d(Z.c_1(\Lcal))>0$ since $\Lcal$ is ample, and we conclude that the left-hand-side of $(\ref{ProjForm})$ is nonzero.  But the right-hand-side of $(\ref{ProjForm})$ vanishes since $\varphi_*(Z)$ is supported on the point $x$.  The contradiction shows that $\varphi$ is finite.

To prove the degree formula, recall that the Euler-Poincar\'e characteristic $\chi(X,\cdot)$ of (tensor) powers of $\Lcal$ satisfies
\begin{equation*}
\chi(X,\Lcal^{\otimes \nu}) = \frac{e(\Lcal)}{\dim(X)!}\nu^{\dim(X)} + \text{ lower order terms}
\end{equation*} 
for some $e(\Lcal)>0$ and all sufficiently large positive integers $\nu$, where the right-hand-side is the Hilbert-Samuel polynomial.  Since $\chi(X,\Lcal^{\otimes d\nu})=\chi(X,\varphi^*\Lcal^{\otimes \nu})=\deg(\varphi)\chi(X,\Lcal^{\otimes \nu})$, comparing leading terms we deduce that
\begin{equation*}
\frac{e(\Lcal)}{\dim(X)!}(d\nu)^{\dim(X)}=\frac{\deg(\varphi)e(\Lcal)}{\dim(X)!}\nu^{\dim(X)}
\end{equation*}
and it follows that $\deg(\varphi)=d^{\dim(X)}$.

{\bf (b)}  Let $Y$ be an irreducible component of the closed subvariety $\Fix(\varphi)$ of $X$.  Note that $Y$ is closed and that $\varphi(Y)=Y$, so $(Y,\varphi|_{Y},\iota^*\Lcal)$ is a polarized algebraic dynamical system, where $\iota:Y\to X$ is the inclusion morphism.  Moreover $\deg(\varphi|_{Y})=d^{\dim(Y)}$ by part {\bf (a)}.  But since $\varphi$ restricted to $Y$ is the identity, we have $\deg(\varphi|_Y)=1$.  It follows that $\dim(Y)=0$ and so the set $\Fix(\varphi)$ is finite.  Since $\Per_n(\varphi)=\Fix(\varphi^n)$, replacing $\varphi$ with $\varphi^n$ we deduce that $\Per_n(\varphi)$ is finite.  Finally, since $\varphi$ is a finite morphism and $\Preper_{n,m}(\varphi)=\varphi^{-m}(\Per_{n}(\varphi))$, we conclude that $\Preper_{n,m}(\varphi)$ is finite. 

{\bf (c)}  This follows from Fakhruddin \cite{Fakhruddin}, who used a result of Hrushovski (\cite{Hrushovski} Thm. 1.1) in model theory to show that $\Per(\varphi)$ is Zariski-dense in $X$.  As remarked in \cite{Fakhruddin}, the larger set $\Preper(\varphi)$ can be shown to be Zariski-dense using the same argument, but without using Hrushovski's theorem.  When $k=\CC$, the Zariski-density of $\Per(\varphi)$ can be proved using results of Briend-Duval \cite{BriendDuval}.
\end{proof}

\end{subsection}

\begin{subsection}{Endomorphisms of $\PP^N_k$}\label{EndoSect}

Let $\varphi:\PP^N_k\to\PP^N_k$ be a surjective morphism; thus $\varphi^*\Ocal(1)\simeq\Ocal(d)$ for some integer $d\geq1$.  If $d=1$ then $\varphi$ is an automorphism, and thus $\deg(\varphi)=1$.  On the other hand, if $d\geq2$ then the triple $(\PP^N_k,\varphi,\Ocal(1))$ is a polarized algebraic dynamical system in the sense of $\S$~\ref{AlgDynSys}, and $\deg(\varphi)=d^N$ by Proposition~\ref{DegreeFormula} {\bf (a)}.

Concretely, choose homogeneous coordinates $\x=(x_0,x_1,\dots,x_N)$ on $\PP^N_{k}$, and let
\begin{equation}\label{HomogMap}
\Phi:k^{N+1}\to k^{N+1} \hskip1cm \Phi(\x) = (\Phi_0(\x),\Phi_1(\x),\dots,\Phi_N(\x))
\end{equation}
be a map defined by $N+1$ homogeneous forms $\Phi_n(\x)\in k[\x]$ of degree $d$.  We say $\Phi$ is {\em nonsingular} if $\Phi(\x)\neq\0$ for all nonzero $\x\in \bar{k}^{N+1}$; in this case $\Phi$ determines a morphism $\varphi:\PP^N_k\to\PP^N_k$ with $\varphi^*\Ocal(1)\simeq\Ocal(d)$.  We call $\Phi$ a {\em model} for $\varphi$ with respect to $\x$, and we will sometimes write this map as $\Phi(\x)$ to indicate the dependence on the choice of coordinates $\x$ on $\PP^N_k$.  Any surjective morphism $\varphi:\PP^N_k\to\PP^N_k$ has such a model $\Phi(\x)$ with respect to $\x$, and if $\Psi(\x)$ and $\Phi(\x)$ are two models for $\varphi$ with respect to $\x$ then $\Psi(\x)=c\Phi(\x)$ for some nonzero constant $c\in k$.

If $\y=(y_0,y_1,\dots,y_N)$ is another choice of coordinates on $\PP^N_k$, then $\Gamma(\x)=\y$ for some $\Gamma\in\GL_{N+1}(k)$.  If $\Psi(\y)$ is a model for a morphism $\varphi$ with respect to the coordinates $\y$, then $\Gamma^{-1}\circ\Psi\circ\Gamma(\x)$ is a model for $\varphi$ with respect to $\x$.

\end{subsection}

\begin{subsection}{Review of the resultant}\label{ResSect}

Fix integers $N\geq1$ and $d\geq1$.  Let $\Phi:k^{N+1}\to k^{N+1}$ be a map defined as in $(\ref{HomogMap})$ by $N+1$ homogeneous forms $\Phi_n(\x)\in k[\x]$ of common degree $d\geq1$ in the variables $\x=(x_0,x_1,\dots,x_N)$.  The resultant $\Res(\Phi)$ of the map is a certain homogeneous integral polynomial in the coefficients of the forms $\Phi_n$; for the definition see \cite{vanderWaerden} $\S$~82 or \cite{Jouanolou}.  For example, when $d=1$ we may view $\Phi$ as an $(N+1)\times(N+1)$ matrix, and $\Res(\Phi)=\det(\Phi)$.  The following proposition states the most basic property of the resultant.

\begin{prop}\label{ResProp}
$\Res(\Phi)=0$ if and only if $\Phi(\x)=\0$ for some nonzero $\x\in \kbar^{N+1}$.
\end{prop}
\begin{proof}
See \cite{vanderWaerden}, $\S$82.
\end{proof}

\end{subsection}

\begin{subsection}{Non-archimedean fields and reduction}\label{NonArchFields}

Throughout this paper $\KK$ denotes a field which is endowed with a nontrivial, non-archimedean absolute value $|\cdot|$.  We denote by $\KK^\circ  = \{\alpha\in \KK\,\mid\,|\alpha|\leq1\}$ the valuation ring of $\KK$, by $\KK^{\circ\circ} = \{\alpha\in \KK\,\mid\,|\alpha|<1\}$ the maximal ideal of $\KK^\circ$, and by $\tilde{\KK}=\KK^\circ/\KK^{\circ\circ}$ the residue field of $\KK$.  For us the most important example occurs when $\KK$ is the function field $K=k(C)$ of a curve $C$ over an algebraically closed constant field $k$, along with the absolute value $|\cdot|_v=e^{-\ord_v(\cdot)}$ associated to a (closed) point $v\in C$.  In this case $\KK^\circ$ coincides with the ring $\Ocal_v$ of regular functions at $v$, and the residue field $\tilde{\KK}$ is isomorphic to the constant field $k$ via the evaluation map $\Ocal_v\to k$.

Let $N\geq1$ be an integer, and let $\x=(x_0,x_1,\dots,x_N)$ denote $N+1$ variables in $\KK$.  Define a norm $\|\cdot\|$ on $\KK^{N+1}$ by $\|\x\|=\max\{|x_0|, |x_1|, \dots, |x_N|\}$, and denote by $B(0,1)=\{\x\in\KK^{N+1}\mid \|\x\|\leq1\}$ the unit ball in $\KK^{N+1}$.  Given a map $F:\KK^{N+1}\to\KK^M$ defined by $M$ polynomials $F_m(\x)\in\KK[\x]$, denote by $H(F)$ the maximum absolute value of the coefficients of $F$; thus $H(F)\leq1$ if and only if $F$ has coefficients in the valuation ring $\KK^\circ$.  

\begin{prop}\label{GLHomogeneous}
Let $\KK$ be an algebraically closed non-archimedean field, and let $F:\KK^{N+1}\to\KK$ be a map defined by a polynomial $F(\x)\in\KK[\x]$.  Then $H(F) = \max\{|F(\x)|\,\mid\, \x\in B(0,1)\}$.
\end{prop}
\begin{proof}
By normalizing $F$ we may assume without loss of generality that $H(F)=1$, and then plainly $|F(\x)|\leq1$ for all $\x\in B(0,1)$ by the ultrametric inequality.  Since $F$ has coefficients in $\KK^\circ$ it reduces to a polynomial $\tilde{F}(\x)\in\tilde{\KK}[\x]$ over the residue field $\tilde{\KK}$.  Since $H(F)=1$, the reduced polynomial $\tilde{F}(\x)$ is nonzero and therefore is nonvanishing on a nonempty Zariski-open subset of $\tilde{\KK}^{N+1}$ (note that $\tilde{\KK}$ is algebraically closed).  Select some $\tilde{\x}_0\in \tilde{\KK}^{N+1}$ such that $\tilde{F}(\tilde{\x}_0)\neq0$, and let $\x_0\in B(0,1)$ be a point which reduces to $\tilde{\x}_0$.  Thus $|F(\x_0)|=1$.
\end{proof}

Let $\Phi:\KK^{N+1}\to \KK^{N+1}$ be a homogeneous map of degree $d\geq1$.  Note that by Proposition~\ref{ResProp}, the map $\Phi$ is nonsingular if and only if $\Res(\Phi)\neq0$.  We say the map $\Phi$ has {\it nonsingular reduction} over $\KK$ if $\Phi$ is defined over $\KK^\circ$ and if the induced map $\tilde{\Phi}:\tilde{\KK}^{N+1}\to \tilde{\KK}^{N+1}$ over the residue field $\tilde{\KK}$ is nonsingular.  By Proposition~\ref{ResProp}, the map $\Phi$ has nonsingular reduction if and only if $\Phi$ has coefficients in $\KK^\circ$ and $|\Res(\Phi)|=1$.  

\begin{lem}\label{HMProp}
Let $\Phi:\KK^{N+1}\to \KK^{N+1}$ be a nonsingular homogeneous polynomial map of degree $d\geq1$.  Then there exist positive constants $C_1,C_2$, depending on $\Phi$, such that $C_1\|\x\|^{d} \leq \|\Phi(\x)\| \leq C_2\|\x\|^{d}$ for all $\x\in \KK^{N+1}$.  If $\Phi$ has coefficients in $\KK^\circ$ then we may take $C_1=|\Res(\Phi)|$ and $C_2=1$.  In particular, if $\Phi$ has nonsingular reduction then $\|\Phi(\x)\| =\|\x\|^{d}$ for all $\x\in\KK^{N+1}$.
\end{lem}
\begin{proof}
The upper bound follows immediately from the ultrametric inequality, and the lower bound follows from basic properties of the resultant; see Prop. 8 of Kawaguchi-Silverman \cite{KawaguchiSilverman}.
\end{proof}

Let $\varphi:\PP^N_\KK\to\PP^N_\KK$ be a morphism of degree at least one.  We say $\varphi$ has {\em good reduction} over $\KK$ if there exists a choice of coordinates $\x=(x_0,x_1,\dots,x_N)$ on $\PP^N_{\KK}$ such that $\varphi$ extends to an endomorphism of the associated integral model $\PP^N_{\KK^\circ}$ of $\PP^N_{\KK}$.  Equivalently, $\varphi$ has good reduction over $\KK$ if there exists a choice of coordinates $\x=(x_0,x_1,\dots,x_N)$ on $\PP^N_{\KK}$, and a model $\Phi(\x)$ for $\varphi$ with respect to $\x$, such that $\Phi(\x)$ has nonsingular reduction as defined above.  Such a model determines a reduced morphism $\tilde{\varphi}:\PP^N_{\tilde{\KK}}\to\PP^N_{\tilde{\KK}}$ over the residue field $\tilde{\KK}$.  

\begin{lem}\label{GammaLemma}
Let $\varphi:\PP^N_\KK\to\PP^N_\KK$ be a morphism of degree at least two.  Let $\Phi(\x)$ and $\Psi(\y)$ be models for $\varphi$ with respect to the coordinates $\x$ and $\y$ on $\PP^N_{\KK}$ respectively, where $\Gamma(\x)=\y$ and $\Phi(\x)=\Gamma^{-1}\circ\Psi\circ\Gamma(\x)$ for some $\Gamma\in\GL_{N+1}(\KK)$.  If both $\Phi(\x)$ and $\Psi(\y)$ have nonsingular reduction, then $\Gamma\in\GL_{N+1}(\KK^\circ)$.
\end{lem}
\begin{proof}
Replacing $\KK$ with $\bar{\KK}$ (and extending the absolute value $|\cdot|$ to $\bar{\KK}$), we may assume without loss of generality that $\KK$ is algebraically closed.  Note that by Lemma~\ref{HMProp} we have $\|\Phi(\x)\|=\|\x\|^d$ for all $\x\in\KK^{N+1}$, since $\Phi(\x)$ has nonsingular reduction, and likewise for $\Psi(\y)$.  By Proposition~\ref{GLHomogeneous} we may select a point $\x_0\in B(0,1)$ where the maximum $H(\Gamma)=\max\{\|\Gamma(\x)\|\,\mid\, \x\in B(0,1)\}$ is achieved.  Therefore
\begin{equation*}
H(\Gamma)^d = \|\Gamma(\x_0)\|^d = \|\Psi(\Gamma(\x_0))\| = \|\Gamma(\Phi(\x_0))\|\leq H(\Gamma),
\end{equation*}
the last inequality following from Proposition~\ref{GLHomogeneous} and the fact that $\Phi(\x_0)\in B(0,1)$.  Since $d\geq2$, we conclude that $H(\Gamma)\leq1$, which means that $\Gamma$ has coefficients in $\KK^\circ$.  By symmetry $\Gamma^{-1}$ has coefficients in $\KK^\circ$ as well, and therefore $\Gamma\in\GL_{N+1}(\KK^\circ)$.
\end{proof}

\end{subsection}

\begin{subsection}{Extending $K$}\label{ExtendSect}

Let $K=k(C)$ and $\varphi:\PP^N_K\to\PP^N_K$ be as in the statement of Theorem~\ref{MainTheorem}.  It is evident from their definitions that the properties of isotriviality and everywhere potential good reduction are invariant under replacing the function field $K=k(C)$ with an extension $K'=k'(C')$ of $K$, where $k'/k$ is an extension of the constant field and $C'\to C$ is a finite map.  Therefore, during the proof of Theorem~\ref{MainTheorem} we may replace $K$ with such an extension $K'$ at any time with no loss of generality.
\end{subsection}

\begin{subsection}{The ``only if'' direction of Theorem~\ref{MainTheorem}}\label{OnlyIfSect}

Let $K=k(C)$ and let $\varphi:\PP^N_K\to\PP^N_K$ be a morphism with $\deg(\varphi)\geq1$.  Assuming that $\varphi$ is isotrivial, it is easy to see that it must have potential good reduction at each place $v\in M_K$.  Extending $K$ if necessary we may assume there exist coordinates $\x=(x_0,x_1,\dots,x_N)$ on $\PP^N_{K}$ such that $\varphi$ has a model $\Phi(\x)$ with respect to $\x$ with coefficients in $k$; thus $\Res(\Phi)$ is a nonzero element of $k$.  Given any place $v$ of $K$, note that $\Phi(\x)$ is defined over $\Ocal_{v}$, since $k\subset\Ocal_{v}$, and $\Res(\Phi)\in k^\times\subset \Ocal_{v}^\times$.  It follows that $\varphi$ has good reduction at $v$. 

\end{subsection}

\begin{subsection}{Automorphisms of $\PP^N_K$}\label{DegOneCaseSect}

Let $K=k(C)$ and let $\varphi:\PP^N_K\to\PP^N_K$ be an automorphism; thus $\varphi^*\Ocal(1)\simeq\Ocal(1)$ and $\deg(\varphi)=1$.  Note that $(\PP^N_K,\varphi,\Ocal(1))$ does not qualify as a polarized algebraic dynamical system in the nomenclature of $\S$~\ref{AlgDynSys}, since it fails the requirement that $d\geq2$.  However, we will show in this section that Theorem~\ref{MainTheorem} still holds in this case.  In place of dynamical tools the proof uses only basic facts from linear algebra.

\begin{proof}[Proof of Theorem~\ref{MainTheorem} for automorphisms]
In view of $\S$~\ref{OnlyIfSect} it suffices to prove the ``if'' part of the theorem.  Let $\varphi:\PP^N_K\to\PP^N_K$ be an automorphism with potential good reduction at all places $v\in M_K$.  Let $\Phi(\x)$ be a model for $\varphi$ with respect to a choice of coordinates $\x=(x_0,x_1,\dots,x_N)$ on $\PP^N_{K}$.  Thus each $\Phi_n(\x)\in K[\x]$ is a linear form, and we may view $\Phi:K^{N+1}\to K^{N+1}$ as a nonsingular $(N+1)\times(N+1)$ matrix over $K$.  Extending $K$ if necessary we may assume without loss of generality that $K$ contains an $(N+1)$-th root of $\det(\Phi)$, and re-normalizing $\Phi$ we may further assume that $\det(\Phi)=1$.  Again extending $K$ if necessary we may assume that $K$ contains all of the eigenvalues of $\Phi$.  Finally, by changing coordinates we may assume that $\Phi$ is in Jordan canonical form.  We are going to show that the eigenvalues of $\Phi$ are in the constant field $k$ of $K$, showing that $\Phi$ is defined over $k$, and completing the proof that $\varphi$ is isotrivial.  

Let $v\in M_K$ be a place of $K$ and let $K'/K$ be a finite extension such that $\varphi$ has good reduction at a place $v'$ of $K'$ over $v$.  It follows that there exists a model $\Psi(\y)$ for $\varphi$, with respect to some choice of coordinates $\y=(y_0,y_1,\dots,y_N)$ on $\PP^N_{K'}$, such that $\Psi$ has coefficients in $\Ocal_{v'}$ and $\det(\Psi)\in\Ocal_{v'}^\times$.  Let $\Gamma\in\GL_{N+1}(K')$ be the change-of-coordinate matrix satisfying $\Gamma(\x)=\y$.  Thus $\Theta(\x):=\Gamma^{-1}\circ\Psi\circ\Gamma(\x)$ is another model for $\varphi$ with respect to the coordinates $\x$, whereby $\Phi(\x)=c\Theta(\x)$ for some nonzero $c\in K'$.  We have $1=\det(\Phi)=c^{N+1}\det(\Theta)=c^{N+1}\det(\Psi)$, and since $\det(\Psi)\in\Ocal_{v'}^\times$ we conclude that $c\in\Ocal_{v'}^\times$ as well.  Letting $P_\Phi(T)=\det(TI-\Phi)$ denote the characteristic polynomial of $\Phi$, and similarly for $\Psi$ and $\Theta$, we have
\begin{equation*}
P_\Phi(T) = P_{c\Theta}(T) = c^{N+1}P_\Theta(T/c) = c^{N+1}P_\Psi(T/c).
\end{equation*}
Since $P_\Psi(T)\in\Ocal_{v'}[T]$ we deduce that $P_\Phi(T)\in\Ocal_{v'}[T]$ as well, and so in fact $P_\Phi(T)\in\Ocal_{v}[T]$ since $P_\Phi(T)$ is defined over the smaller field $K$.  Since $P_\Phi(T)\in\Ocal_{v}[T]$ is monic and splits over $K$ we conclude that the eigenvalues of $\Phi$ are in $\Ocal_v$.  As $v\in M_K$ was arbitrary, the eigenvalues of $\Phi$ are in $\Ocal_v$ at all $v\in M_K$, and so they must be in the constant field $k$ as desired.
\end{proof}

\end{subsection}

\end{section}


\begin{section}{The Geometric Proof of Theorem~\ref{MainTheorem}}\label{GeomSect}

\begin{subsection}{Overview}  

Throughout  this section we let $k$ denote an algebraically closed field, and we fix integers $N\geq1$ and $d\geq2$.  In this section we will study the space $\End_{N,d}(k)$ of endomorphisms $\varphi$ of $\PP^N_k$ with $\varphi^*\Ocal(1)\simeq\Ocal(d)$.  Generalizing a result of Silverman \cite{SilvermanSpace}, we will show in $\S$~\ref{GeomQuot} that the quotient $\Mcal_{N,d}(k)=\End_{N,d}(k)/\PGL_{N+1}(k)$ of this space by the automorphism group of $\PP^N_k$ is an affine $k$-variety.  In $\S$~\ref{GeomProofSect} we will give the geometric proof of Theorem~\ref{MainTheorem}.  

\end{subsection}

\begin{subsection}{The space of endomorphisms}\label{space}

Let $\Sym^d(k^{N+1})^{N+1}$ be the space of homogeneous maps $k^{N+1}\to k^{N+1}$ of degree $d$.  Explicitly, an element of this space is given by an $(N+1)$-tuple $\Phi(\x) = (\Phi_0(\x),\Phi_1(\x),\dots,\Phi_N(\x))$, where each $\Phi_n(\x)\in k[\x]$ is a homogeneous form of degree $d$ in the variables $\x=(x_0,x_1,\dots,x_N)$.

Given a choice of homogeneous coordinates $\x=(x_0,x_1,\dots,x_N)$ on $\PP^N_k$, recall from $\S$~\ref{EndoSect} that each morphism $\varphi:\PP^N_k\rightarrow \PP^N_k$ with $\varphi^*\Ocal(1)\simeq\Ocal(d)$ has a model $\Phi(\x)\in \Sym^d(k^{N+1})^{N+1}$, which is unique up to scaling by a constant $c\in k^\times$, and that moreover $\Res(\Phi)\neq0$.  Thus $\varphi$ corresponds to a unique point in the projective space $\PP((\Sym^d(k^{N+1}))^{N+1})$.  Note that since $\Res(\Phi)$ is itself a homogeneous form in the coefficients of $\Phi$, the condition $\Res(\Phi)=0$ defines a closed hypersurface $\Res_{N,d}(k)$ in $\PP((\Sym^d(k^{N+1}))^{N+1})$ (see \cite{GKZ}, $\S$~3.3).

In view of these remarks and Proposition~\ref{DegreeFormula} {\bf (a)}, we define the {\em space of endomorphisms of} $\PP^N_k$ {\em of degree} $d^N$ by
\begin{equation*}
\End_{N,d}(k) := \PP((\Sym^d(k^{N+1}))^{N+1})\setminus \Res_{N,d}(k).
\end{equation*}
Thus $\End_{N,d}(k)$ is an affine open subvariety of $\PP((\Sym^d(k^{N+1}))^{N+1})$.

Note that the correspondence between morphisms $\varphi:\PP^N_k\rightarrow \PP^N_k$ of degree $d^N$ and points in $\End_{N,d}(k)$ depends on our initial choice of coordinates $\x$.  Moreover, changing coordinates on $\PP^N_k$ corresponds to conjugating $\varphi$ by an element $\gamma$ of $\PGL_{N+1}(k)$, and so we are led to consider the action of $\PGL_{N+1}(k)$ on $\End_{N,d}(k)$ by $(\gamma,\varphi)\mapsto \gamma^{-1}\varphi\gamma$.  Therefore, the ``correct'' coordinate-independent space parametrizing morphisms $\varphi:\PP^N_k\rightarrow \PP^N_k$ of degree $d^N$ is the quotient 
\begin{equation}\label{SetQuotient}
\Mcal_{N,d}(k) := \End_{N,d}(k)/\PGL_{N+1}(k)
\end{equation}
of this action.  While this quotient can be defined set-theoretically, it is not guaranteed that $\Mcal_{N,d}(k)$ is a variety over $k$, or that the fibers of the quotient map $\End_{N,d}(k)\to\Mcal_{N,d}(k)$ are closed.  In the next section we will show that the quotient $\Mcal_{N,d}(k)$ does in fact naturally carry the structure of a $k$-variety, and that the associated quotient map is a morphism.

\end{subsection}

\begin{subsection}{Existence of a geometric quotient}\label{GeomQuot}  We recall a basic definition from geometric invariant theory.  Let $\alpha:G\times X\to X$ be an action of an algebraic group $G$ over $k$ on a $k$-variety $X$.  A pair $(Y,\pi)$ consisting of a $k$-variety $Y$ and a morphism $\pi:X\to Y$ is called a {\em geometric quotient} of $X$ by the action of $G$ if it satisfies the following properties:
\begin{quote}
\textbf{(i)}  The diagram
\begin{equation*}
\begin{CD}
G\times X  @> \alpha >>  X \\ 
@V p_2 VV                                    @VV \pi V \\ 
X @> \pi >>  Y
\end{CD}
\end{equation*}
commutes.  \\
\textbf{(ii)}  $\pi$ is surjective, and the image of $(\alpha,p_2):G\times X\rightarrow X\times X$ is $X\times_Y X$. \\
\textbf{(iii)}  A subset $U\subset Y$ is open if and only if $\pi^{-1}(U)$ is open in $X$. \\
\textbf{(iv)}  The fundemental sheaf $\Ocal_Y$ is the subsheaf of $\pi_*(\Ocal_X)$ consisting of invariant functions.
\end{quote}
We refer to \cite{Dolgachev} and \cite{MFK} for additional definitions and background material.

The purpose of this section is to show that the action of $\PGL_{N+1}(k)$ on $\End_{N,d}(k)$ has a geometric quotient.  We begin by recording several preliminary results.  Recall that the stablizer of a point $\varphi\in\End_{N,d}(k)$ by the action of $\PGL_{N+1}(k)$ is the subgroup $\St(\varphi)=\{\gamma\in\PGL_{N+1}(k)\mid\gamma^{-1}\varphi\gamma=\varphi\}$.  Given a subset $S$ of $\PP^N_k$, we say $S$ is in general position if every nonempty finite subset $T$ of $S$ with $|T|\leq N+1$ is linearly independent.  The following lemma is well-known.

\begin{lem}\label{PGLProp}
Let $S$ be a subset of $\PP^N_k$ in general position with $|S|=N+2$ points.  If $\gamma\in\PGL_{N+1}(k)$ is an automorphism of $\PP^N_k$ and $\gamma(P)=P$ for all $P\in S$, then $\gamma$ is the identity automorphism.
\end{lem}

\begin{prop}\label{stabilizerfinite}
Let $\varphi\in\End_{N,d}(k)$.  The stablizer $\St(\varphi)$ of $\varphi$ by the action of $\PGL_{N+1}(k)$ is finite.
\end{prop}

\begin{proof}
By Proposition~\ref{DegreeFormula} {\bf (c)}, the set $\Preper(\varphi)$ of $\varphi$-preperiodic points in $\PP^N_k$ is Zariski-dense.  It follows that there exists a set $S=\{P_0,\dots,P_{N+1}\}$ of $N+2$ preperiodic points in general position.  Explicitly, if $r\leq N-1$, and if $P_0,\dots, P_r$ are $r+1$ linearly independent preperiodic points, we choose a point $P_{r+1}$ in the projective $N$-space which is preperiodic and not in the linear space of dimension $r$ generated by the points $P_0,\dots, P_r$; this is possible because $\Preper(\varphi)$ is Zariski-dense.  Applying this for $r=0,1,\dots,N-1$, we obtain $N+1$ linearly independent preperiodic points $P_0,\dots,P_N$.  Again using the Zariski-density of $\Preper(\varphi)$, we let $P_{N+1}$ be a preperiodic point which is not on any of the $N+1$ hyperplanes generated by $N$-point subsets of $\{P_0,\dots,P_N\}$.

Each preperiodic point $P_i$ lies in $\Preper_{n_i,m_i}(\varphi)$ for some integers $n_i\geq1$ and $m_i\geq1$.  Note that given $\gamma\in\St(\varphi)$, we have $\gamma^{-1}\varphi\gamma=\varphi$, which implies $\gamma^{-1}\varphi^r\gamma=\varphi^r$ and $\gamma\varphi^r=\varphi^r\gamma$ for any positive integer $r$.  If $P$ is in $\Preper_{n_i,m_i}(\varphi)$ then
\begin{equation*}
\varphi^{m_i}\gamma(P)=\gamma\varphi^{m_i}(P)=\gamma\varphi^{n_i+m_i}(P)=\varphi^{n_i+m_i}\gamma(P),
\end{equation*}
so $\gamma(P)\in \Preper_{n_i,m_i}(\varphi)$ as well.  Thus $\St(\varphi)$ acts on each finite set $\Preper_{n_i,m_i}(\varphi)$ and we obtain a group homomorphism
\begin{equation}\label{GroupHom}
\St(\varphi)\to\prod^{N+1}_{i=0}\Perm\left(\Preper_{n_i,m_i}(\varphi)\right),
\end{equation}
where $\Perm\left(\Preper_{n_i,m_i}(\varphi)\right)$ denotes the group of permutations of the set $\Preper_{n_i,m_i}(\varphi)$.  If $\gamma$ is in the kernel of the map $(\ref{GroupHom})$, then in particular it fixes each point in the set $S=\{P_0,P_1,\dots,P_{N+1}\}$, whereby $\gamma$ is the identity automorphism by Lemma~\ref{PGLProp}.  Thus the map $(\ref{GroupHom})$ is injective, and since each set $\Preper_{n_i,m_i}(\varphi)$ is finite, it follows that $\St(\varphi)$ is finite.
\end{proof}

\begin{cor}\label{actionclosed}
The action of $\PGL_{N+1}(k)$ on $\End_{N,d}(k)$ is closed.
\end{cor}

\begin{proof}
By an argument on p. 10 of \cite{MFK}, if for each $\varphi\in \End_{N,d}(k)$ there exists an open neighborhood $U$ of $\varphi$ where the dimension of the stabilizer $\St(\psi)$ is constant for all $\psi\in U$, then the action by $\PGL_{N+1}(k)$ is closed.  Since $\St(\varphi)$ is zero-dimensional for all $\varphi\in\End_{N,d}(k)$ by Proposition~\ref{stabilizerfinite}, the action by $\PGL_{N+1}(k)$ is closed.
\end{proof}

\begin{prop}\label{GQexistance}
A geometric quotient of $\End_{N,d}(k)$ by $\PGL_{N+1}(k)$ exists, and moreover it is affine.
\end{prop}

\begin{proof}
The proof is just an application of Amplification 1.3 of \cite{MFK}.  Recall $\End_{N,d}(k)$ is affine and $\PGL_{N+1}(k)$ is reductive.  Therefore an affine geometric quotient exists if and only if the action of $\PGL_{N+1}(k)$ is closed, which is the case by Corollary \ref{actionclosed}.
\end{proof}

We let $(\Mcal_{N,d},\pi)$ denote the geometric quotient of $\End_{N,d}(k)$ by the action of $\PGL_{N+1}(k)$.  Note that the set of points $\Mcal_{N,d}(k)$ on this quotient coincides with the set-theoretic quotient defined in $(\ref{SetQuotient})$.  We remark that the case $N=1$ of Proposition~\ref{GQexistance} follows immediately from Silverman \cite{SilvermanSpace}, Thm. 1.1.

\end{subsection}

\begin{subsection}{The geometric proof of Theorem~\ref{MainTheorem}}\label{GeomProofSect}  In view of $\S$~\ref{OnlyIfSect} and $\S$~\ref{DegOneCaseSect}, it suffices to consider the ``if'' direction of the statement for morphisms $\varphi:\PP^N_K\to\PP^N_K$ with $\varphi^*\Ocal(1)\simeq\Ocal(d)$ for $d\geq2$.  Moreover, by the remarks of $\S$~\ref{ExtendSect} we may assume without loss of generality that the constant field $k$ of $K=k(C)$ is algebraically closed.  

Suppose that $\varphi$ has potential good reduction at all places $v$ of $K$.  Since $\varphi$ fails to have good reduction at only finitely many places of $K$, by extending $K$ we may assume without loss of generality that $\varphi$ has good reduction at all places $v$ of $K$.

Let $v$ be a place of $K$; in other words $v\in C$ is a (closed) point.  Recall that the surjective map $\Ocal_v\to k$ given by evaluation at $v$ induces a canonical isomorphism between the residue field of $\Ocal_v$ and the constant field $k$.  Since $\varphi$ has good reduction at $v$ there exists a model $\Phi_v(\x)$ for $\varphi$ with respect to a choice of coordinates $\x=(x_0,x_1,\dots,x_N)$ on $\PP^N_{K}$, such that the coefficients of $\Phi_v(\x)$ are in $\Ocal_v$ and $\Res(\Phi_v)\in\Ocal_v^\times$.  Reduction modulo the maximal ideal of $\Ocal_v$ defines a map $\tilde{\Phi}_v:k^{N+1}\to k^{N+1}$ and an associated morphism $\tilde{\varphi}_v:\PP^N_k\to\PP^N_k$.

Moreover, the model $\Phi_v(\x)$ has nonsingular reduction at all but finitely many points $u\in C$.  Therefore, we can find an affine open neighborhood $U_v=\Spec(A_v)$ of $v\in C$ such that $\Phi_v(\x)$ has coefficients in $A_v\subset K$, and such that $\Phi_v(\x)$ has nonsingular reduction at all points $u\in U_v$.  As at the point $v$, reduction at each point $u\in U_v$ defines a morphism $\tilde{\varphi}_u:\PP^N_k\to\PP^N_k$.  We obtain a morphism $U_v\rightarrow\End_{N,d}$ defined by $u\mapsto\tilde{\varphi}_u$.  The sets $\{U_v\}_{v\in C}$ define an affine open cover of $C$, and since the curve $C$ is quasi-compact we can find a finite subcover $\{U_i\}$.  Thus for each open set $U_i$ in our finite cover of $C$ we have morphisms
\begin{equation}\label{UtoM}
U_i\rightarrow\End_{N,d}\xrightarrow{\pi}\Mcal_{N,d}.
\end{equation}

Let $U_{ij}=U_i\cap U_j$ be an intersection between two of the open sets in the finite cover, and let $\Phi_i(\x)$ and $\Phi_j(\y)$ be the models for $\varphi$ with respect to the neighborhoods $U_i$ and $U_j$ as described above, respectively.  Let $\Gamma\in\GL_{N+1}(K)$ denote the change of coordinate element satisfying $\Gamma(\x)=\y$.  Thus $\Phi_i(\x)=c\Gamma^{-1}\circ\Phi_j\circ\Gamma(\x)$ for some $c\in K^\times$.  Extending $K$ if necessary we may assume there exists some $a\in K$ such that $a^{d-1}=c$.  Letting $\Gamma'=a\Gamma$, we have $\Phi_i(\x)=(\Gamma')^{-1}\circ\Phi_j\circ\Gamma'(\x)$.  Therefore, replacing $\Gamma$ with $\Gamma'$ and replacing the coordinates $\y$ with $\y'=\Gamma'(\x)=a\y$, we may assume without loss of generality that $\Phi_i(\x)=\Gamma^{-1}\circ\Phi_j\circ\Gamma(\x)$.  

Given an arbitrary point $u\in U_{ij}=U_i\cap U_j$, both $\Phi_i(\x)$ and $\Phi_j(\y)$ have nonsingular reduction at $u$, which means that $\Gamma\in \GL_{N+1}(\Ocal_u)$ by Lemma~\ref{GammaLemma}.  Denote by $\tilde{\Gamma}_u\in\GL_{N+1}(k)$ the reduction of $\Gamma$ at $u$ and by $\tilde{\gamma}_u\in\PGL_{N+1}(k)$ the associated automorphism of $\PP_k^N$; thus $\tilde{\Phi}_{u,i}(\x)=\tilde{\Gamma}_u^{-1}\circ\tilde{\Phi}_{u,j}\circ\tilde{\Gamma}_u(\x)$.  Therefore, letting $\tilde{\varphi}_{u,i}\in \End_{N,d}$ denote the endomorphism of $\PP_k^N$ obtained by reduction of the model $\Phi_i(\x)$ at $u$, and likewise defining $\tilde{\varphi}_{u,j}\in \End_{N,d}$ using the model $\Phi_j(\y)$ at $u$, we deduce that $\tilde{\varphi}_{u,i}=\tilde{\gamma}_u^{-1}\circ\tilde{\varphi}_{u,j}\circ\tilde{\gamma}_u$.  We have shown that the image of $u\in U_{ij}$ in $\End_{N,d}$ is well-defined up to $\PGL_{N+1}(k)$-conjugation; that is, it is contained in a unique fiber of the quotient map $\pi$.  We obtain a morphism $C\rightarrow\Mcal_{N,d}$ by the inclusion $U_i\hookrightarrow C$ and (\ref{UtoM}).  

The quotient $\Mcal_{N,d}$ is affine by Proposition \ref{GQexistance} and $C$ is complete.  Hence the image of $C\rightarrow\Mcal_{N,d}$ is a point.  By \cite{Dolgachev} Corollary 6.1, the fiber of this point in $\Mcal_{N,d}$ contains a unique closed $\PGL_{N+1}(k)$-conjugacy class.  It follows that $\varphi$ coincides with the base extension $\psi_K:\PP^N_K\to\PP^N_K$ for some $\psi$ in this class, which means that $\varphi$ is isotrivial.
\end{subsection}

\end{section}


\begin{section}{The Analytic Proof of Theorem~\ref{MainTheorem}}

\begin{subsection}{Homogeneous transfinite diameter}\label{HomTDSect}

Let $\KK$ be a non-archimedean field as discussed in $\S$~\ref{NonArchFields}, and let $E$ be a bounded, infinite subset of $\KK^{N+1}$.  In this section we will define the homogeneous transfinite diameter $d_\infty(E)$, a nonnegative real number which in a certain sense measures the size of $E$.  In the two dimensional case ($N=1$) this variation on the classical notion of transfinite diameter was introduced and studied by Baker-Rumely \cite{BakerRumely}.  When $N\geq2$ our definition of $d_\infty(E)$ is new.

To define $d_\infty(E)$, let $M\geq N+1$ be an integer, and let $\Scal_M(E)$ denote the set of subsets $S$ of $E$ with exactly $|S|=M$ elements.  Given $S\in \Scal_M(E)$, enumerate by $S_1,\dots S_{J_M}$ the subsets of $S$ with exactly $|S_j|=N+1$ elements; thus $J_M={\binom{M}{N+1}}$.  Define
\begin{equation}\label{Vandermonde}
\Delta(S)=\prod_{1\leq j\leq J_{M}}\det(S_j),
\end{equation}
where $\det(S_j)$ is the determinant of the $(N+1)\times(N+1)$ matrix whose column-vectors are the elements of $S_j$.  Thus each $\det(S_j)$ and the product $(\ref{Vandermonde})$ are defined only up to sign.  Define the $M$-diameter of $E$ by
\begin{equation}\label{MDiam}
d_M(E) = \sup_{S\in\Scal_M(E)}|\Delta(S)|^{1/J_M}.
\end{equation}
A standard argument shows that the sequence $d_M(E)$ is monotone decreasing as $M\to+\infty$, and therefore the limit 
\begin{equation}\label{TransDiam}
d_\infty(E)=\lim_{M\to+\infty}d_M(E)
\end{equation}
exists; see Proposition \ref{TDProp1} {\bf (a)} below.  The number $d_\infty(E)$ is called the (homogeneous) transfinite diameter of $E$.  The following proposition summarizes some basic properties of the transfinite diameter.  Given a set $E\subset \KK^{N+1}$, we say that $E$ is an {\em ellipsoid} if it is of the form $E=\Gamma(B(0,1))$ for some $\Gamma\in\GL_{N+1}(\KK)$.

\begin{prop}\label{TDProp1}
Let $\KK$ be a non-archimedean field, and let $E$ be a bounded, infinite subset of $\KK^{N+1}$.
\begin{quote}
{\bf (a)}  The sequence $d_M(E)$ is monotone decreasing, and thus the limit $(\ref{TransDiam})$ exists. \\
{\bf (b)}  If $\Gamma\in \GL_{N+1}(\KK)$, then $d_\infty(\Gamma(E))=|\det(\Gamma)|d_\infty(E)$.
\end{quote}
For the remainder of this proposition, assume that $\KK$ is algebraically closed.
\begin{quote}
{\bf (c)}  If $E$ contains the unit ball $B(0,1)$, then $d_\infty(E)\geq1$; moreover $d_\infty(B(0,1))=1$. \\
{\bf (d)}  If $E$ is an ellipsoid such that $B(0,1)\subseteq E$ and $d_\infty(E)=1$, then $E=B(0,1)$. \\
{\bf (e)}  If $E$ is an ellipsoid which contains the standard unit basis elements $\e_0=(1,0,\dots,0)$, $\e_1=(0,1,\dots,0)$, $\dots$, $\e_N=(0,\dots,0,1)$, then $B(0,1)\subseteq E$. 
\end{quote}
\end{prop}
\begin{proof}
{\bf (a)}  The following is a variation on the standard argument for the existence of the transfinite diameter; it generalizes the proof given for $N=1$ in \cite{BakerRumely}, Lemma 3.10.  Fix $M\geq N+1$ and $\epsilon>0$.  By the definition $(\ref{MDiam})$ we may choose a set $S=\{\x(1), \dots, \x(M+1)\}$ of $M+1$ elements in $E$ such that $|\Delta(S)|\geq(d_{M+1}(E)-\epsilon)^{J_{M+1}}$.  For each $1\leq m \leq M+1$ denote by $T_m=S\setminus\{\x(m)\}$; thus $|T_m|=M$ and $|\Delta(T_m)|\leq d_M(E)^{J_M}$ by $(\ref{MDiam})$.  Observe that
\begin{equation*}
\prod_{1\leq m\leq M+1} |\Delta(T_m)| = |\Delta(S)|^{M-N},
\end{equation*}
since the left-hand-side is the product of $|\det(S_j)|$ over each $(N+1)$-element subset $S_j$ of $S$ exactly $(M+1)-(N+1)=M-N$ times, which is precisely the same as right-hand-side.  Thus
\begin{equation*}
d_{M+1}(E)-\epsilon \leq |\Delta(S)|^{1/J_{M+1}} \leq d_M(E)^{J_M(M+1)/J_{M+1}(M-N)} = d_M(E).
\end{equation*}
Since $\epsilon>0$ is arbitrary, we conclude that $d_{M+1}(E)\leq d_M(E)$.

{\bf (b)}  Since $\Delta(\Gamma(S))=\pm\det(\Gamma)^{J_M}\Delta(S)$, we have $d_M(\Gamma(E))=|\det(\Gamma)|d_M(E)$ for all $M$, and the claim follows.

{\bf (c)}  Suppose that $B(0,1)\subseteq E$.  Let $\{\tilde{\x}(m)\}_{m=1}^{\infty}$ be an infinite sequence of points in $\tilde{\KK}^{N+1}$ such that any $N+1$ terms of the sequence are linearly independent over $\tilde{\KK}$.  [To see that such a sequence exists, let $\tilde{\x}(1), \dots, \tilde{\x}(N+1)$ be any basis for $\tilde{\KK}^{N+1}$; we define the rest of the sequence by induction.  Suppose that $m\geq N+1$ and that the first $m$ terms $\tilde{\x}(1), \dots, \tilde{\x}(m)$ of the sequence have been constructed with the desired linear-independence property.  Each choice of $N$ elements in the set $\{\tilde{\x}(1), \dots, \tilde{\x}(m)\}$ spans a hyperplane in $\tilde{\KK}^{N+1}$; since $\KK$ (and therefore $\tilde{\KK}$), is algebraically closed, we let $\tilde{\x}(m+1)$ be any element in the complement of the union of these hyperplanes.  Clearly any $N+1$ elements in the set $\{\tilde{\x}(1), \dots, \tilde{\x}(m+1)\}$ are linearly independent, and by induction on $m$ the sequence exists as claimed.]

For each $m\geq1$ let $\x(m)\in B(0,1)$ be a point reducing to $\tilde{\x}(m)$, and fix $M\geq N+1$.  Let $S=\{\x(1), \dots, \x(M)\}\in\Scal_M(E)$, and note that by the $\tilde{\KK}$-linear-independence property of the sequence $\{\tilde{\x}(m)\}_{m=1}^{\infty}$ we have $|\Delta(S)|=1$.  It follows from $(\ref{MDiam})$ that $d_M(E)\geq1$ for all $M\geq N+1$, and the inequality $d_\infty(E)\geq1$ follows from $(\ref{TransDiam})$.

If $E=B(0,1)$, then the opposite inequality $d_\infty(E)\leq1$ follows at once from $(\ref{MDiam})$, $(\ref{TransDiam})$, and the ultrametric inequality.

{\bf (d)}  Since $E$ is an ellipsoid we have $E=\Gamma(B(0,1))$ for some $\GL_{N+1}(\KK)$, and since $B(0,1)\subseteq E$, we conclude that $\Gamma^{-1}(B(0,1))\subseteq B(0,1)$.  Thus $\Gamma^{-1}$ maps $B(0,1)$ into itself, and it follows from Proposition~\ref{GLHomogeneous} that $\Gamma^{-1}$ has coefficients in $\KK^\circ$.  Also, $1=d_\infty(E)=d_\infty(\Gamma(B(0,1)))=|\det(\Gamma)|d_\infty(B(0,1))=|\det(\Gamma)|$ by part {\bf (b)} of this Proposition.  We conclude that $\Gamma^{-1}$, and therefore also $\Gamma$, is an element of $\GL_{N+1}(\KK^\circ)$, from which it follows that $E=\Gamma(B(0,1))=B(0,1)$.

{\bf (e)}  Again $E=\Gamma(B(0,1))$ for some $\Gamma\in\GL_{N+1}(\KK)$, and thus $\Gamma^{-1}(\e_n)\in B(0,1)$ for all $n$.  Thus given an arbitrary $\x\in B(0,1)$ we have
\begin{equation*}
\x = \sum_{n=0}^{N}x_n\e_n = \Gamma(\sum_{n=0}^{N}x_n\Gamma^{-1}(\e_n)) \in \Gamma(B(0,1)) = E;
\end{equation*}
here the containment follows from the ultrametric inequality and the fact that $x_n\in\KK^\circ$ and $\Gamma^{-1}(\e_n) \in B(0,1)$.  Thus $B(0,1)\subseteq E$.
\end{proof}

\end{subsection}

\begin{subsection}{Homogeneous local height functions and filled Julia sets}\label{HMSection}

We will now define two dynamical objects associated to each homogeneous map $\Phi:\KK^{N+1}\to \KK^{N+1}$ of degree $d\geq2$.  The {\it homogeneous local height function} $\Hhat_\Phi:\KK^{N+1}\setminus\{\0\}\to\RR$ associated to $\Phi$ is defined by
\begin{equation}\label{HLHF}
\Hhat_\Phi(\x) = \lim_{\ell\to+\infty}\frac{1}{d^\ell}\log\|\Phi^\ell(\x)\|,
\end{equation}
where $\Phi^\ell=\Phi\circ\dots\circ\Phi$ denotes the map $\Phi$ composed with itself $\ell$ times.  It is easy to show using Lemma~\ref{HMProp} that $\Hhat_\Phi$ defines a continuous real-valued function on $\KK^{N+1}\setminus\{\0\}$, and we may extend the definition of $\Hhat_\Phi$ by declaring that $\Hhat_\Phi(\0)=-\infty$.  These functions were defined in the case $N=1$ by Baker-Rumely \cite{BakerRumely}, and were generalized and further studied by Kawaguchi-Silverman in \cite{KawaguchiSilverman}, \cite{KawaguchiSilverman2}.  The {\it homogeneous filled Julia set} associated to $\Phi$ is the set 
\begin{equation}\label{HFJS}
F_\Phi = \{\x\in \KK^{N+1} \,\mid\, \sup_{\ell\geq1}\|\Phi^\ell(\x)\|<+\infty \}
\end{equation}
of points whose forward iterates remain bounded.  We will summarize the basic properties of $\Hhat_\Phi$ and $F_\Phi$ in the following proposition.

\begin{prop}\label{HFJSProp}
Let $\KK$ be a non-archimedean field, let $\Phi:\KK^{N+1}\to \KK^{N+1}$ be a nonsingular homogeneous map of degree $d\geq2$, let $c\in \KK^\times$, and let $\Gamma\in\GL_{N+1}(\KK)$.  Then the following identities hold:
\begin{quote}
{\bf (a)}  $\Hhat_\Phi(c\x) = \Hhat_\Phi(\x) + \log|c|$; \\
{\bf (b)}  $\Hhat_{c\Phi}(\x) = \Hhat_\Phi(\x) + \frac{1}{d-1}\log|c|$; \\
{\bf (c)}  $\Hhat_{\Gamma^{-1}\circ\Phi\circ\Gamma}(\x) =\Hhat_\Phi(\Gamma(\x))$; \\
{\bf (d)}  $F_\Phi = \{\x\in \KK^{N+1} \,\mid\, \Hhat_\Phi(\x)\leq0\}$; \\
{\bf (e)}  $F_{c\Phi} = c^{-1/(d-1)}F_\Phi$, assuming $c^{-1/(d-1)}\in\KK$; \\
{\bf (f)}  $F_{\Gamma^{-1}\circ\Phi\circ\Gamma} =\Gamma^{-1}(F_\Phi)$.
\end{quote}
\end{prop}
\begin{proof}  {\bf (a)} This follows immediately from the definition $(\ref{HLHF})$.  

{\bf (b)} Note that $(c\Phi)^\ell(\x)=c^{d^{\ell-1}+d^{\ell-2}+\dots d+1}\Phi^\ell(\x)=c^{(d^\ell-1)/(d-1)}\Phi^\ell(\x)$, and therefore
\begin{equation*}
\frac{1}{d^\ell}\log\|(c\Phi)^\ell(\x)\| = \frac{1}{d^\ell}\log\|\Phi^\ell(\x)\| + \frac{d^\ell-1}{d^\ell(d-1)}\log|c|;
\end{equation*}
letting $\ell\to+\infty$ establishes the desired identity.

{\bf (c)}  By Lemma~\ref{HMProp} applied to the map $\Gamma^{-1}:\KK^{N+1}\to \KK^{N+1}$ we have 
\begin{equation*}
\frac{1}{d^\ell}\log\|(\Gamma^{-1}\circ\Phi\circ\Gamma)^\ell(\x)\| = \frac{1}{d^\ell}\log\|\Gamma^{-1}(\Phi^\ell(\y))\|= \frac{1}{d^\ell}(\log\|\Phi^\ell(\y)\|+O(1))
\end{equation*}
where $\y=\Gamma(\x)$ and $O(1)$ denotes a function which is bounded as $\ell\to+\infty$; letting $\ell\to+\infty$ establishes the desired identity.

{\bf (d)}  If $\x\in F_\Phi$ then the iterates $\Phi^\ell(\x)$ are bounded, and $\Hhat_\Phi(\x)\leq0$ follows immediately from the definition $(\ref{HLHF})$.  Conversely, suppose that $\x\not\in F_\Phi$.  Let $T>0$ be a parameter, and by assumption we have $\|\Phi^{\ell_0}(\x)\|>T$ for some $\ell_0$ depending on $T$.  For each $\ell>\ell_0$ we then have
\begin{equation*}
\begin{split}
\|\Phi^{\ell}(\x)\| & \geq C_1^{1+d+d^2+\dots+d^{\ell-\ell_0-1}} \|\Phi^{\ell_0}(\x)\|^{d^{\ell-\ell_0}} \\
	& > C_1^{1+d+d^2+\dots+d^{\ell-\ell_0-1}} T^{d^{\ell-\ell_0}} \\
	& = C_1^{(d^{\ell-\ell_0}-1)/(d-1)} T^{d^{\ell-\ell_0}}
\end{split}
\end{equation*}
by iterating the lower bound in Lemma \ref{HMProp}.  Thus
\begin{equation*}
\begin{split}
\Hhat_\Phi(\x) & = \lim_{\ell\to+\infty}\frac{1}{d^\ell}\log\|\Phi^\ell(\x)\| \\
	& \geq \lim_{\ell\to+\infty}\frac{1}{d^\ell}\log(C_1^{(d^{\ell-\ell_0}-1)/(d-1)} T^{d^{\ell-\ell_0}}) \\
	& = d^{-\ell_0}(d^{-1}\log C_1 + \log T).
\end{split}
\end{equation*}
Selecting any $T>C_1^{1/d}$ we deduce that $\Hhat_\Phi(\x)>0$.

{\bf (e)}  This follows at once from {\bf (d)} along with {\bf (a)} and {\bf (b)}.

{\bf (f)}  This follows at once from {\bf (d)} and {\bf (c)}.
\end{proof}

The following proposition characterizes the property of nonsingular reduction for a homogeneous map in terms of its homogeneous local height function and its homogeneous filled Julia set.

\begin{prop}\label{GoodReductionProp}
Let $\KK$ be an algebraically closed non-archimedean field, and let $\Phi:\KK^{N+1}\to \KK^{N+1}$ be a nonsingular homogeneous map of degree $d\geq2$.  The following four conditions are equivalent:
\begin{quote}
	{\bf (a)}  $\Phi$ has nonsingular reduction; \\
	{\bf (b)}  $\|\Phi(\x)\|=\|\x\|^d$ for all $\x\in \KK^{N+1}$; \\
	{\bf (c)}  $\Hhat_\Phi(\x)=\log\|\x\|$ for all $\x\in \KK^{N+1}$; \\
	{\bf (d)}  $F_\Phi=B(0,1)$.
\end{quote}
\end{prop}

{\em Remarks.}  This proposition generalizes Lemma 3.9 of \cite{BakerRumely} from $N=1$ to arbitrary $N\geq1$.  The equivalence of {\bf (a)} and {\bf (c)} was proved by Kawaguchi-Silverman in \cite{KawaguchiSilverman}, Prop. 14; our proof that {\bf (d)} implies {\bf (a)} borrows a key argument from their paper.

\begin{proof}
If {\bf (a)} holds then {\bf (b)} follows immediately from Lemma~\ref{HMProp}.  That {\bf (b)} implies {\bf (c)} is clear from the definition $(\ref{HLHF})$.  If {\bf (c)} holds, then {\bf (d)} follows immediately from Proposition \ref{HFJSProp} {\bf (d)}.

Finally, suppose that {\bf (d)} holds.  Plainly the filled Julia set always satisfies $\Phi(F_\Phi)\subseteq F_\Phi$, which in this case implies that $\|\Phi(\x)\|\leq1$ for all $\x\in B(0,1)$.  Thus by Proposition~\ref{GLHomogeneous} we have $H(\Phi)\leq1$, which means that $\Phi$ is defined over $\KK^\circ$.  In particular, it follows that $|\Res(\Phi)|\leq1$ by the ultrametric inequality.  

Suppose that $|\Res(\Phi)|<1$.  Then by Proposition~\ref{ResProp} (taking $k$ to be the residue field $\tilde{\KK}=\KK^\circ/\KK^{\circ\circ}$) there exists $\x_0\in B(0,1)$ such that $\x_0\not\equiv0\pmod{\KK^{\circ\circ}}$ but $\Phi(\x_0)\equiv0\pmod{\KK^{\circ\circ}}$.  In particular, we have $\|\x_0\|=1$ and $\|\Phi(\x_0)\|<1$.  By iterating the upper bound in Proposition~\ref{HMProp} we have $\|\Phi^k(\x_0)\| \leq \|\Phi(\x_0)\|^{d^{k-1}}$, and thus
\begin{equation*}
\Hhat_\Phi(\x_0)=\lim_{k\to\infty}\frac{1}{d^k}\log\|\Phi^k(\x_0)\| \leq \frac{1}{d}\log\|\Phi(\x_0)\|<0.
\end{equation*}
Since $\KK$ is algebraically closed we can find a (nonzero) scalar $c\in \KK$ such that $\Hhat_\Phi(\x_0)<\log|c|<0$.  In particular $|c|<1$, whereby $\|c^{-1}\x_0\|>1$ and thus $c^{-1}\x_0\notin B(0,1)$.  On the other hand by Proposition~\ref{HFJSProp} {\bf (a)} we have
\begin{equation*}
\Hhat_\Phi(c^{-1}\x_0)=\Hhat_\Phi(\x_0) - \log|c| <0,
\end{equation*}
which according to Proposition~\ref{HFJSProp} {\bf (d)} implies that $c^{-1}\x_0$ is in the filled Julia set $F_\Phi$.   These two properties of $c^{-1}\x_0$ contradict the assumption {\bf (d)} that $F_\Phi=B(0,1)$.  We conclude that $|\Res(\Phi)|=1$, and we have shown that {\bf (d)} implies {\bf (a)}.
\end{proof}

The following lemma calculates the transfinite diameter of the filled Julia set of certain homogeneous maps.

\begin{lem}\label{TDJulia}
Let $\KK$ be an algebraically closed non-archimedean field, and let $\Phi:\KK^{N+1}\to \KK^{N+1}$ be a homogeneous map of degree $d\geq2$.  Suppose that a conjugate $\Psi=\Gamma^{-1}\circ\Phi\circ\Gamma$ of $\Phi$ by some $\Gamma\in\GL_{N+1}(\KK)$ has nonsingular reduction.  Then
\begin{equation*}
d_\infty(F_\Phi) = |\Res(\Phi)|^{C(N,d)},
\end{equation*}
where $C(N,d)$ is a constant depending only on $N$ and $d$.
\end{lem}

\begin{proof}
We will need to use the following composition law for the resultant: if $\Phi$ and $\Phi'$ are systems of $N+1$ homogeneous forms in $N+1$ variables, of degrees $d$ and $d'$ respectively, then 
\begin{equation}\label{ResCompLaw}
\Res(\Phi\circ\Phi') = \Res(\Phi)^{a}\Res(\Phi')^{b},
\end{equation}
where $a$ and $b$ are constants depending only on $N$, $d$, and $d'$; for a proof of this fact see \cite{ChengMckayWang}, Cor. 5.

Turning now to proof of the lemma, since $\Psi=\Gamma^{-1}\circ\Phi\circ\Gamma$ has nonsingular reduction, we have $F_\Psi=B(0,1)$, and so
\begin{equation*}
1=d_\infty(F_\Psi) =d_\infty(\Gamma^{-1}(F_\Phi)) = |\det(\Gamma)|^{-1}d_\infty(F_\Phi)
\end{equation*}
by Proposition~\ref{TDProp1} {\bf (b)} and {\bf (c)}.  On the other hand
\begin{equation*}
1=|\Res(\Psi)| = |\Res(\Gamma^{-1}\circ\Phi\circ\Gamma)|=|\det(\Gamma)|^{A(N,d)}|\Res(\Phi)|^{B(N,d)},
\end{equation*}
by $(\ref{ResCompLaw})$ and the fact that $\Res(\Gamma)=\det(\Gamma)$, where $A(N,d)$ and $B(N,d)$ are constants depending only on $N$ and $d$.  Thus $d_\infty(F_\Phi) = |\det(\Gamma)| = |\Res(\Phi)|^{-B(N,d)/A(N,d)}$.
\end{proof}

{\em Remarks.}  The hypothesis in Lemma \ref{TDJulia} that some conjugate of $\Phi$ have nonsingular reduction is probably unnecessary, although we do not know a proof for general $\Phi$.  A proof of this identity for general $\Phi$ is given in the case $\KK=\CC_p$ by DeMarco-Rumely \cite{DeMarcoRumely}, although we do not know whether their proof generalizes to the equal-characteristic case.  We also point out that, using an explicit form of the composition law $(\ref{ResCompLaw})$, as in say \cite{ChengMckayWang}, it is possible to give an explicit expression for the exponent $C(N,d)$, although we will not need to do so.

\end{subsection}


\begin{subsection}{The analytic proof of Theorem~\ref{MainTheorem}}\label{AnalyticProof}

Let $K=k(C)$ be a function field as discussed in the introduction, let $M_K$ denote the set of places of $K$, or equivalently, the set of closed points on the curve $C$.  For each $v\in M_K$ denote by $|\cdot|_v$ a non-archimedean absolute value on $K$ associated to $v$, normalized so that the product formula holds in the form $\prod_{v\in M_K}|a|_v =1$ for all nonzero $a\in K$.  Let $\KK_v$ be an algebraically closed non-archimedean field containing $K$ equipped with an absolute value extending $|\cdot|_v$.  Thus $\Ocal_v\subset\KK_v^\circ$.

We now give the analytic proof of Theorem~\ref{MainTheorem}, which is stated in a slightly stronger form here.  This result and its proof generalizes Thm. 1.9 of Baker \cite{Baker} and Prop. 6.1 of Benedetto \cite{Benedetto}.

\begin{thm}\label{IsoThm}
Let $K=k(C)$ be a function field, and let $\varphi:\PP^N_K\to\PP^N_K$ of degree at least two.  The following are equivalent:
\begin{quote}
	{\bf (a)}  $\varphi$ is isotrivial; \\
	{\bf (b)}  $\varphi$ has potential good reduction at all places $v\in M_K$; \\
	{\bf (c)}  $\varphi$ has good reduction over $\KK_v$ for all places $v\in M_K$.
\end{quote}
\end{thm}

\begin{proof}
We proved that {\bf (a)} implies {\bf (b)} in $\S$~\ref{OnlyIfSect}.  It is trivial that {\bf (b)} implies {\bf (c)}, since if $v$ is any place of $K$, if $K'/K$ is any finite extension, and if $v'$ is any place of $K'$ lying over $v$, then there exists an embedding $K'\hookrightarrow\KK_v$ with $\Ocal_{v'}\hookrightarrow\KK_v^\circ$.

Finally, we will show that {\bf (c)} implies {\bf (a)}.  Suppose that $\varphi$ has good reduction over $\KK_v$ at all places $v$ of $K$.  Note that both conditions {\bf (a)} and {\bf (c)} are invariant under replacing $K$ with an extension of $K$ as in $\S$~\ref{ExtendSect}.  By Proposition~\ref{DegreeFormula} {\bf (c)}, the $\varphi$-preperiodic points are Zariski-dense in $\PP^N(\Kbar)$.  Therefore, by extending $K$ if necessary, we may assume that there exist at least $N+1$ linearly independent $\varphi$-preperiodic points in $\PP^N(K)$.  Choose coordinates $\x=(x_0,x_1,\dots,x_N)$ on $\PP^N_K$ such that these $N+1$ $\varphi$-preperiodic points are the points $P_0,P_1,\dots,P_N\in\PP^N(K)$ which lift to the standard basis elements $\e_0=(1,0,\dots,0)$, $\e_1=(0,1,0,\dots,0),\dots$, $\e_N=(0,\dots,0,1)$.

Without loss of generality, we may also assume that there exists a model $\Phi(\x)$ for $\varphi$ with respect to the coordinates $\x$ such that the standard basis elements $\e_0$, $\e_1$, $\dots$, $\e_N\in K^{N+1}$ are $\Phi$-preperiodic.  To see this, let $\Psi:K^{N+1}\to K^{N+1}$ be {\it any} model for $\varphi$ with respect to $\x$, and note that for each $n$, since $P_n$ is $\varphi$-preperiodic we have $\varphi^{i_n}(P_n)=\varphi^{j_n}(P_n)$ for some integers $1\leq i_n<j_n$.  Thus $\Psi^{i_n}(\e_n)=c_n\Psi^{j_n}(\e_n)$ for some nonzero constant $c_n\in K$.  For each $n$ select an element $\alpha_n\in \Kbar$ such that $\alpha_n^{d^{i_n}}c_n=\alpha_n^{d^{j_n}}$; thus 
\begin{equation*}
\begin{split}
\Psi^{i_n}(\alpha_n \e_n) & = \alpha_n^{d^{i_n}}\Psi^{i_n}(\e_n) \\
	&  = \alpha_n^{d^{i_n}}c_n\Psi^{j_n}(\e_n) = \alpha_n^{d^{j_n}}\Psi^{j_n}(\e_n) = \Psi^{j_n}(\alpha_n\e_n).
\end{split}
\end{equation*}
Therefore each $\alpha_n\e_n$ is $\Psi$-preperiodic.  Replace $K$ by a finite extension containing the $\alpha_n$, and let 
$\Phi(\x')=\Gamma^{-1}\circ\Psi\circ\Gamma(\x')$, where $\Gamma\in\GL_{N+1}(K)$ is selected to take $\e_n$ to $\alpha_n\e_n$ for each $n$.  Now the standard basis elements $\e_n$ are $\Phi$-preperiodic.  Replacing the coordinates $\x$ with $\x'=\Gamma^{-1}(\x)$, the above claim is justified.

To summarize, we have an endomorphism $\varphi$ of $\PP^N_K$ which has good reduction over $\KK_v$ for all places $v$ of $K$, a choice of coordinates $\x=(x_0,x_1,\dots,x_N)$ on $\PP^N_K$, and a model $\Phi(\x)$ for $\varphi$ with respect to $\x$ such that the standard basis elements $\e_0, \e_1 , \dots, \e_N$ of $K^{N+1}$ are $\Phi$-preperiodic.  We are going to show that
\begin{equation}\label{FPhiOhat}
F_{\Phi,v} = B_v(0,1) \,\text{ for all }\,\,v\in M_K,
\end{equation}
where $F_{\Phi,v}$ denotes the homogeneous filled Julia set in $\KK_v^{N+1}$ associated to $\Phi(\x)$, and $B_v(0,1)$ denotes the unit ball in $\KK_v^{N+1}$.  Assuming this claim, it follows from Proposition~\ref{GoodReductionProp} that $\Phi(\x)$ has nonsingular reduction at all places $v\in M_K$, which means in particular that the coefficients of $\Phi(\x)$ are in $\KK_v^\circ$ for all $v\in M_K$.  This implies that the coefficients of $\Phi(\x)$ are elements of the constant field $k$ of $K$, since a rational function on $C$ with no poles must be constant.  Therefore $\varphi$ is defined over $k$, whereby it is isotrivial, completing the proof that {\bf (c)} implies {\bf (a)}. 

It now remains only to prove $(\ref{FPhiOhat})$.  Fix a place $v\in M_K$.  Since $\varphi$ has good reduction over $\KK_v$ there exists a choice of coordinates $\y=(y_0:y_1:\dots:y_N)$ on $\PP^N_{\KK_v}$ and a model $\Psi(\y)$ for $\varphi$ with respect to $\y$ such that $\Psi(\y)$ has nonsingular reduction over $\KK_v$; thus $\Psi(\y)$ has coefficients in $\KK_v^\circ$ and $|\Res(\Psi)|_v=1$.  Moreover $F_{\Psi,v}=B_v(0,1)$ by Proposition~\ref{GoodReductionProp}.

Choose $\Gamma\in\GL_{N+1}(\KK_v)$ so that $\y=\Gamma(\x)$.  Thus $\Phi'(\x)=\Gamma^{-1}\circ\Psi\circ\Gamma(\x)$ is another model for $\varphi$ with respect to the coordinates $\x$, so $\Phi'(\x)=c\Phi(\x)$ for some $c\in\KK_v^\times$.  Therefore $\Gamma(c^{-1/(d-1)}F_{\Phi,v})=\Gamma(F_{\Phi',v})=F_{\Psi,v}=B_v(0,1)$ by Proposition~\ref{HFJSProp} {\bf (e)} and {\bf (f)}; thus $F_{\Phi,v}=c^{1/(d-1)}\Gamma^{-1}(B_v(0,1))$.  In particular, $F_{\Phi,v}$ is an ellipsoid, as defined in $\S$\ref{HomTDSect}.  Note also that the standard basis elements $\e_n$ are elements of $F_{\Phi,v}$ since they are $\Phi$-preperiodic.  By Proposition~\ref{TDProp1} {\bf (e)} we conclude that $B_v(0,1)\subseteq F_{\Phi,v}$, which implies by Proposition \ref{TDProp1} {\bf (c)} that  $d_\infty(F_{\Phi,v})\geq1$.

On the other hand, by Lemma \ref{TDJulia} and the product formula we have 
\begin{equation*}
\prod_{v\in M_K}d_\infty(F_{\Phi,v}) = \prod_{v\in M_K}|\Res(\Phi)|_v^{C(N,d)}=1,
\end{equation*}  
and since we have already shown that $d_\infty(F_{\Phi,v})\geq1$ for all $v\in M_K$, we must actually have $d_\infty(F_{\Phi,v})=1$ for all $v\in M_K$.  Since each $F_{\Phi,v}$ is an ellipsoid containing $B_v(0,1)$, we deduce from Proposition \ref{TDProp1} {\bf (d)} that $F_{\Phi,v} = B_v(0,1)$ for all $v\in M_K$.  Thus we have proved $(\ref{FPhiOhat})$, which completes the proof that {\bf (c)} implies {\bf (a)}.
\end{proof}

\end{subsection}

\end{section}


\begin{section}{Two Applications}\label{ApplicationsSect}

\begin{subsection}{Endomorphisms with an isotrivial iterate}

The following corollary of Theorem~\ref{MainTheorem} states that an endomorphism is isotrivial if and only if it has an isotrivial iterate.

\begin{cor}\label{IterateCor}
Let $K=k(C)$ be a function field, let $\varphi:\PP^N_K\to\PP^N_K$ be a morphism of degree at least two, and let $r\geq1$ be an integer.  Then $\varphi$ is isotrivial if and only if $\varphi^r$ is isotrivial.
\end{cor}

\begin{proof}[Proof of Corollary~\ref{IterateCor}]
By the equivalence of {\bf (a)} and {\bf (c)} in Theorem~\ref{IsoThm}, it suffices to show that given any place $v\in M_K$, $\varphi$ has good reduction over $\KK_v$ if and only if $\varphi^r$ has good reduction over $\KK_v$.  The ``only if'' direction of this statement is trivial.  To show the ``if'' direction, suppose that $\varphi^r$ has good reduction over $\KK_v$.  Thus there exist coordinates $\x=(x_0,x_1,\dots,x_N)$ on $\PP^N_{\KK_v}$, and a model $\Psi(\x)$ for $\varphi^r$ with respect to $\x$ such that $\Psi(\x)$ has nonsingular reduction over $\KK_v$; thus $F_{\Psi,v}=B_v(0,1)$ by Proposition~\ref{GoodReductionProp}.  Let $\Phi(\x)$ be a model for $\varphi$ with respect to the same coordinates $\x$; thus $\Phi^r(\x)$ is a model for $\varphi^r$, so $\Phi^r(\x)=c\Psi(\x)$ for some $c\in \KK_v^\times$.  It follows from Proposition~\ref{HFJSProp} {\bf (e)} that
\begin{equation*}
F_{\Phi,v}=F_{\Phi^r,v}=F_{c\Psi,v}=c^{-1/(d-1)}F_{\Psi,v}=c^{-1/(d-1)}B_v(0,1).
\end{equation*}
Letting $\Phi'=c^{-1}\Phi$, we have $F_{\Phi'_v} = c^{1/(d-1)}F_{\Phi,v}=B(0,1)$ by Proposition~\ref{HFJSProp} {\bf (e)}, whereby $\Phi'(\x)$ is a model for $\varphi$ with nonsingular reduction by Proposition~\ref{GoodReductionProp}.  Therefore $\varphi$ has good reduction, as desired.
\end{proof}

\end{subsection}

\begin{subsection}{A dynamical criterion for decomposability of locally free coherent sheaves}\label{BundleSect}

Let $C$ be a complete nonsingular curve over an algebraically closed field $k$.  Let $N\geq1$ be an integer, let $\Ecal$ be a locally free coherent sheaf of rank $N+1$ on $C$, and denote by $\pi:\PP(\Ecal)\to C$ the associated projective bundle.  The following corollary of Theorem~\ref{MainTheorem} states that, after possibly replacing $C$ with a base extension $p:C'\to C$ and replacing $\Ecal$ with $\Ecal'=p^*\Ecal$, the sheaf $\Ecal$ decomposes as a direct sum of $N+1$ copies of the same invertible sheaf on $C$ if and only if there exists an endomorphism of $\PP(\Ecal)$ of degree at least two.  A similar result was obtained by Amerik \cite{Amerik} in the case $k=\CC$.

\begin{cor}
Let $\Ecal$ be a locally free coherent sheaf of rank $N+1$ on a complete nonsingular curve $C$ over an algebraically closed field $k$.  Then the following two conditions are equivalent:
\begin{quote}
{\bf (a)}  There exists a base extension $p:C'\to C$ and an endomorphism $\varphi:\PP(\Ecal')\to\PP(\Ecal')$ of degree at least two; \\
{\bf (b)}  There exists a base extension $p:C'\to C$ and an invertible sheaf $\Lcal$ on $C'$ such that $\Ecal'\simeq\Lcal\oplus\dots\oplus\Lcal$.
\end{quote}
Moreover, if {\bf (a)} and {\bf (b)} hold then the two extensions $C'$ can be chosen to coincide, and $\PP(\Ecal')\simeq\PP^N_k\times C'$ with $\varphi=\varphi_0\times\Id_{C'}$, where $\varphi_0:\PP^N_k\to\PP^N_k$ is a morphism and $\Id_{C'}:C'\to C'$ is the identity.
\end{cor}
\begin{proof}
Both conditions {\bf (a)} and {\bf (b)}  are invariant under replacing $C$ with a finite extension $p:C'\to C$ (and replacing $\Ecal$ with $\Ecal'=p^*\Ecal$), and therefore we may do this at any time with no loss of generality.  Moreover, $\PP(\Ecal)\simeq\PP(\Ecal\otimes\Bcal)$ for any invertible sheaf $\Bcal$ on $C$, so we may also replace $\Ecal$ with $\Ecal\otimes\Bcal$ at any time with no loss of generality.  We identify the generic fiber of $\pi:\PP(\Ecal)\to C$ with $\PP^N_K$, where $K=k(C)$ denotes the function field of $C$, and given a morphism $\varphi:\PP(\Ecal)\to\PP(\Ecal)$ we denote by $\varphi_K:\PP^N_K\to\PP^N_K$ the restriction of $\varphi$ to $\PP^N_K$.  

First suppose that {\bf (b)} holds.  Replacing $C$ with a suitable extension $p:C'\to C$ we may assume that $\Ecal\simeq\Lcal\oplus\dots\oplus\Lcal$ for some invertible sheaf $\Lcal$ on $C$.  Moreover, replacing $\Ecal$ with $\Ecal\otimes\Lcal^\vee$ we may assume that $\Ecal\simeq\Ocal_C\oplus\dots\oplus\Ocal_C$ is isomorphic to the trivial vector bundle.  Thus $\PP(\Ecal)\simeq\PP^N_k\times C$, and any endomorphism $\varphi_0:\PP^N_k\to\PP^N_k$ of degree at least two induces such an endomorphism $\varphi=\varphi_0\times\Id_{C}$ of $\PP(\Ecal)\simeq\PP^N_k\times C$, completing the proof that {\bf (b)} implies {\bf (a)}.

Conversely, suppose that {\bf (a)} holds.  We have $\varphi^*\Ocal_{\PP(\Ecal)}(1)\simeq\Ocal_{\PP(\Ecal)}(d)\otimes\pi^*\Acal$ for some $d\geq2$ and $\Acal\in\Pic(C)$.  Replacing $C$ with a suitable extension $p:C'\to C$ we may assume there exists $\Bcal\in\Pic(C)$ such that $\Bcal^{\otimes (1-d)}\simeq\Acal$, and it follows that $\varphi^*\Ocal_{\PP(\Ecal\otimes\Bcal)}(1)\simeq\Ocal_{\PP(\Ecal\otimes\Bcal)}(d)$.   Replacing $\Ecal$ with $\Ecal\otimes\Bcal$, we may assume without loss of generality that $\varphi^*\Ocal_{\PP(\Ecal)}(1)\simeq\Ocal_{\PP(\Ecal)}(d)$.  

Since $\PP(\Ecal)$ is locally isomorphic to $\PP^N_k\times U$ for open sets $U\subset C$, the morphism $\varphi_K$ has everywhere good reduction.  Therefore by Theorem~\ref{MainTheorem} it is isotrivial, which means that, after replacing $C$ with a suitable extension $p:C'\to C$ if necessary, $\varphi_K$ is induced by an endomorphism $\varphi_0:\PP^N_k\to\PP^N_k$.  In particular, there exist coordinates $\x=(x_0,\dots,x_N)$ on $\PP^N_K$ and a model $\Phi(\x)$ for $\varphi_K:\PP^N_K\to\PP^N_K$ with coefficients in the constant field $k$.

Given a point $P\in\PP^N(K)\subset\PP^N_K\subset\PP(\Ecal)$ and a closed point $v\in C$, the valuative criterion for properness (\cite{Hartshorne} Thm. II.4.7) determines a unique point $s_P(v)\in\pi^{-1}(v)$ specializing $P$.  This defines a section $s_P:C\to\PP(\Ecal)$ of $\pi$, along with a surjective morphism $\Ecal\to s_P^*\Ocal_{\PP(\Ecal)}(1)$ of sheaves on $C$.  Moreover, if $P$ is a $\varphi$-preperiodic point then, after perhaps replacing $C$ with an extension $p:C'\to C$, we have 
\begin{equation}\label{PrePerSheafIdentity}
s_P^*\Ocal_{\PP(\Ecal)}(1)\simeq\Ocal_C.
\end{equation}
To see this note that $s_{\varphi(P)}^*\Ocal_{\PP(\Ecal)}(1)\simeq s_{P}^*\varphi^*\Ocal_{\PP(\Ecal)}(1)\simeq s_{P}^*\Ocal_{\PP(\Ecal)}(d)$.  Thus if $P$ is $\varphi$-preperiodic with $\varphi^{n+m}(P)=\varphi^{m}(P)$ for $n\geq1$ and $m\geq0$, then $s_{P}^*\Ocal_{\PP(\Ecal)}(d^{n+m})\simeq s_{P}^*\Ocal_{\PP(\Ecal)}(d^m)$.  This implies that $s_{P}^*\Ocal_{\PP(\Ecal)}(d^{n+m}-d^{m})\simeq \Ocal_C$, which means that $s_{P}^*\Ocal_{\PP(\Ecal)}(1)$ is a torsion element of $\Pic(C)$.  After replacing $C$ with a suitable extension $p:C'\to C$ we deduce $(\ref{PrePerSheafIdentity})$ as desired.

Now let $P_0,P_1,\dots,P_N\in\PP^N(K)\subset\PP^N_K$ be a linearly independent set of $k$-rational $\varphi$-preperiodic points; such a set exists since $\varphi_K$ is defined over $k$ and since Proposition~\ref{DegreeFormula} {\bf (c)} ensures that the preperiodic points are Zariski-dense in $\PP^N(k)$.  Extending $C$ if necessary we may assume that $(\ref{PrePerSheafIdentity})$ holds for each $P\in\{P_0,P_1,\dots,P_N\}$, and we obtain a morphism
\begin{equation}\label{MorSheaves}
\Ecal\to\bigoplus_{j=0}^{N}s_{P_j}^*\Ocal_{\PP(\Ecal)}(1)\simeq\Ocal_C\oplus\dots\oplus\Ocal_C
\end{equation}
of sheaves on $C$.  In fact we are going to show that $(\ref{MorSheaves})$ is an isomorphism; for this it suffices to show that the set $\{s_{P_0}(v), \dots, s_{P_N}(v)\}$ is linearly independent on each closed fiber $\pi^{-1}(v)\simeq\PP^N_k$ of $\pi:\PP(\Ecal)\to C$.  

Given a point $v\in C$, there exists a neighborhood $U\subset C$ of $v$ such that $\pi^{-1}(U)\simeq\PP^N_k\times U$ and a model $\Psi(\y)$ for $\varphi_K$ which coincides with the morphism $\varphi:\PP(\Ecal)\to\PP(\Ecal)$ when restricted to $\pi^{-1}(U)$.  In particular, $\Psi(\y)$ has nonsingular reduction at $v$.  Let $\Gamma\in\GL_{N+1}(K)$ denote the change of coordinate element satisfying $\Gamma(\x)=\y$.  Thus $\Phi(\x)=c\Gamma^{-1}\circ\Psi\circ\Gamma(\x)$ for some $c\in K^\times$.  Extending the curve $C$ if necessary we may assume there exists some $a\in K$ such that $a^{d-1}=c$.  Letting $\Gamma'=a\Gamma$, we have $\Phi(\x)=(\Gamma')^{-1}\circ\Psi\circ\Gamma'(\x)$.  Therefore, replacing $\Gamma$ with $\Gamma'$ and replacing the coordinates $\y$ with $\y'=\Gamma'(\x)=a\y$, we may assume without loss of generality that $\Phi(\x)=\Gamma^{-1}\circ\Psi\circ\Gamma(\x)$.  By Lemma~\ref{GammaLemma}, since both $\Phi(\x)$ and $\Psi(\y)$ have nonsingular reduction, $\Gamma$ must in fact be an element of $\GL_{N+1}(\Ocal_v)$, which means it reduces to an automorphism $\gamma_v:\PP^N_k\to\PP^N_k$ over the residue field $k$ at $v$.  Moreover $\gamma_v(P_j)=s_{P_j}(v)$ for all $0\leq j\leq N$, and since the set $\{P_0,\dots,P_N\}$ is linearly independent in $\PP^N(k)$ it follows that the set $\{s_{P_0}(v), \dots, s_{P_N}(v)\}$ is linearly independent in $\PP^N(k)$ as well.  Since this holds for all $v\in C$ we deduce that $(\ref{MorSheaves})$ is an isomorphism.  Thus $\Ecal\simeq\Ocal_C\oplus\dots\oplus\Ocal_C$, completing the proof that {\bf (a)} implies {\bf (b)}.  Since $\Ecal$ is the trivial sheaf we have $\PP(\Ecal)\simeq\PP^N_k\times C$ with $\varphi=\varphi_0\times\Id_C$.
\end{proof}

\end{subsection}

\end{section}


\medskip


\begin{thebibliography}{1}

\bibitem{Amerik}
E. Amerik.
\newblock On endomorphisms of projective bundles.
\newblock {\em Man. Math.} 111 (1), (2003), 17-28.

\bibitem{Baker}
M. Baker.
\newblock A finiteness theorem for canonical heights attached to rational maps over function fields. 
\newblock preprint (2005), to appear in {\em J. Reine Angew. Math.}

\bibitem{BakerRumely}
M. Baker and R. Rumely.
\newblock Equidistribution of small points, rational dynamics, and potential theory.
\newblock {\em Ann. Inst. Fourier (Grenoble).}  56, no. 3 (2006), 625--688.

\bibitem{Benedetto}
R. Benedetto.
\newblock Heights and preperiodic points of polynomials over function fields.
\newblock {\em Int. Math. Research Notices}, 62, (2005), 3855-3866.

\bibitem{BriendDuval} 
J.-Y. Briend and J. Duval.
\newblock Exposants de Liapounoff et distribution des points périodiques d'un endomorphisme de $\CC\PP^n$.
\newblock {\em Acta Math.}, 182 (1999), No. 2, 143-157.

\bibitem{Cartan} 
H. Cartan, et. al.
\newblock Fonctions automorphes.
\newblock {\em S\'eminaire Henri Cartan}, 10 (2), (1958).

\bibitem{ChambertLoir} 
A. Chambert-Loir.
\newblock Th\'eor\`emes d'\'equidistribution pour les syst\`emes dynamiques d'origine arithm\'etique.
\newblock Preprint (2006), to appear in {\em Panoramas et synth\`ses}.

\bibitem{ChatzidakisHrushovski} 
Z. Chatzidakis and E. Hrushovski.
\newblock Difference fields and descent in algebraic dynamics - I.
\newblock Preprint (2007).

\bibitem{ChengMckayWang}
C. Cheng, J. Mckay, and S. Wang.
\newblock A chain rule for multivariable resultants.
\newblock {\em Proc. Am. Math. Soc.}, 123(4), (1995), 1037-1047.

\bibitem{DeMarcoRumely}
L. DeMarco and R. Rumely.
\newblock Transfinite diameter and the resultant.
\newblock To appear in {\em J. Reine Angew. Math.}

\bibitem{Dolgachev} 
I. Dolgachev.
\newblock Lectures on Invariant Theory.
\newblock Cambridge University Press, 2003.
London Mathematic Society Lecture Notes Series, No. 296.

\bibitem{Fakhruddin}
N. Fakhruddin.
\newblock Questions on self maps of algebraic varieties. 
\newblock {\em J. Ramanujan Math. Soc.} 111(2) (2003), 109-122. 

\bibitem{GKZ} 
I. M. Gel'fand, M. M. Kapranov, and A. V. Zelevinsky.
\newblock Discriminants, Resultants and Multidimensional Determinants
\newblock Birkhauser, Boston 1994.

\bibitem{Hartshorne} 
R. Hartshorne.
\newblock Algebraic Geometry.
\newblock Springer-Verlag,  New York, 1977. Graduate Texts in Mathematics, No. 52.

\bibitem{Hrushovski} 
E. Hrushovski.
\newblock The elementary theory of the Frobenius automorphisms.
\newblock Preprint (2004).

\bibitem{Jouanolou} 
J. P. Jouanolou.
\newblock Le formalisme du r\'esultant.
\newblock {\em Advances in Mathematics.} 90, (1991), 117-263.

\bibitem{KawaguchiSilverman}
S. Kawaguchi and J. H. Silverman.
\newblock Dynamics of projective morphisms having identical canonical heights.
\newblock {\em Proc. London Math. Soc.} 95(2), (2007), 519-544.

\bibitem{KawaguchiSilverman2}
S. Kawaguchi and J. H. Silverman.
\newblock Nonarchimedean Green functions and dynamics on projective space.
\newblock preprint (2007), to appear in {\em Math. Zeit.}

\bibitem{Kodaira}
K. Kodaira.
\newblock A certain type of irregular algebraic surfaces.
\newblock {\em Journal d'Analyse Math\'ematique} 19(1) (1967) 207-215.

\bibitem{MFK} 
D. Mumford, J. Fogarty and F. Kirwan.
\newblock Geometric Invariant Theory, Third Enlarged Edition.
\newblock Springer-Verlag, New York, 1994.
Ergebnisse der Mathematik und ihrer Grenzgebiete, No. 34.

\bibitem{SilvermanSpace} 
J. Silverman.
\newblock The space of rational maps on $\PP^1$.
\newblock {\em Duke Math. J.}, Volume 94, Number 1 (1998), 41-77.

\bibitem{SilvermanDynamics} 
J. Silverman.
\newblock The Arithmetic of Dynamical Systems.
\newblock Springer-Verlag,  New York, 2007. Graduate Texts in Mathematics, No. 241.

\bibitem{vanderWaerden}
B. van der Waerden.
\newblock {\em Modern Algebra. Vol. II}.
\newblock Frederick Ungar Publishing Co., New York, 1949.

\bibitem{Zhang}
S. Zhang.
\newblock Distributions  in Algebraic Dynamics,
\newblock {\em Survey in Differential Geometry} 10, 381-430, International Press, 2006

\end{thebibliography}
\end{document}